\documentclass[12pt]{amsart}
\usepackage{amsmath,amsthm,amsfonts,amscd,amssymb,eucal,latexsym,mathrsfs,enumerate}

\numberwithin{equation}{section}

\theoremstyle{plain}
\newtheorem{lemma}{Lemma}[section]
\newtheorem{proposition}[lemma]{Proposition}
\newtheorem{theorem}[lemma]{Theorem}
\newtheorem{corollary}[lemma]{Corollary}

\theoremstyle{definition}
\newtheorem{definition}[lemma]{Definition}
\newtheorem{example}[lemma]{Example}
\newtheorem{question}[lemma]{Question}

\theoremstyle{remark}

\newcommand{\GN}{\mathcal{GN}}
\newcommand{\N}{\mathcal{N}}
\newcommand{\vnotimes}{\overline{\otimes}}
\newcommand{\Aut}{\mathrm{Aut}}
\newcommand{\Fix}{\mathrm{Fix}\,}
\newcommand{\ad}{\mathrm{Ad}}

\begin{document}

\title[Groupoid normalisers]{Groupoid normalisers of tensor products: Infinite von Neumann algebras}

\author{Junsheng Fang}
\address{\hskip-\parindent
Junsheng Fang, Department of Mathematics, Texas A\&{}M University, College Station, Texas, 77843, U.S.A.}
\email{jfang@math.tamu.edu}
\author{Roger R.~Smith}
\address{\hskip-\parindent
Roger R.~Smith, Department of Mathematics, Texas A{\&}M University,
College Station, Texas 77843, U.S.A.}
\email{rsmith@math.tamu.edu}
\author{Stuart White}
\address{\hskip-\parindent
Stuart White, Department of Mathematics, University of Glasgow, 
University Gardens, Glasgow Q12 8QW, U.K.}
\email{s.white@maths.gla.ac.uk}

\date{8 March, 2010}
\maketitle

\begin{abstract}
The groupoid normalisers of a unital inclusion $B\subseteq M$ of von Neumann algebras consist of the set $\GN_M(B)$ of partial isometries $v\in M$ with $vBv^*\subseteq B$ and $v^*Bv\subseteq B$. Given two unital inclusions $B_i\subseteq M_i$ of von Neumann algebras, we examine groupoid normalisers for the tensor product inclusion $B_1\ \vnotimes\ B_2\subseteq M_1\ \vnotimes\ M_2$ establishing the formula
$$
\GN_{M_1\,\vnotimes\,M_2}(B_1\ \vnotimes\ B_2)''=\GN_{M_1}(B_1)''\ \vnotimes\ \GN_{M_2}(B_2)''
$$
when one inclusion has a discrete relative commutant $B_1'\cap M_1$ equal to the centre of $B_1$ (no assumption is made on the second inclusion). This result also holds when one inclusion is a generator masa in a free group factor.  We also examine when a unitary $u\in M_1\ \vnotimes\ M_2$ normalising a tensor product $B_1\ \vnotimes\ B_2$ of irreducible subfactors factorises as $w(v_1\otimes v_2)$ (for some unitary $w\in B_1\ \vnotimes\ B_2$ and normalisers $v_i\in\N_{M_i}(B_i)$). We obtain a positive result when one of the $M_i$ is finite or both of the $B_i$ are infinite. For the remaining case, we characterise the II$_1$ factors $B_1$ for which such factorisations always occur (for all $M_1, B_2$ and $M_2$) as those with a trivial fundamental group.
\end{abstract}

\section{Introduction}

This paper is concerned with the behaviour of structrural properties of inclusions of  von Neumann algebras obtained from tensor products. This important construction has a rich history,   and we mention two particularly significant results.  Given two inclusions $B_i\subseteq M_i$ ($i=1,2$) of von Neumann algebras, Tomita's commutation theorem (see \cite{Tak.Tom}) determines the relative commutant of $B_1\ \vnotimes\ B_2$ in $M_1\ \vnotimes\ M_2$ showing that
\begin{equation}\label{eqq1.1}
(B_1\ \vnotimes\ B_2)'\cap(M_1\ \vnotimes\ M_2)=(B_1'\cap M_1)\ \vnotimes\ (B_2'\cap M_2).
\end{equation}
In a different direction, Ge and Kadison \cite{Ge.Kadison} showed that when  $M_1$ and $M_2$ are factors,  any von Neumann algebra lying between $M_1\ \overline{\otimes} \ 1$ and $M_1\  \vnotimes\ M_2$ must have the form $M_1\ \vnotimes\ B_2$ for a von Neumann algebra $B_2\subseteq M_2$.
Equation (\ref{eqq1.1}) can be interpreted as saying that the two operations of taking tensor products and passing to relative commutants commute. In \cite{Saw.Groupoid}, Wiggins and the authors examined 
commutation questions of this form for the operations of taking tensor products and passing to the algebra generated by the groupoid normalisers in the context of inclusions of finite von Neumann algebras.  This paper examines these questions for general inclusions of von Neumann algebras.

The normalisers of an inclusion $B\subseteq M$ of von Neumann algebras consist of those unitaries $u\in M$ with $uBu^*=B$. These elements form a group denoted  $\N_M(B)$ and the von Neumann algebra $\N_M(B)''$ they generate is the \emph{normalising algebra} of $B$ in $M$.  Normalisers were first studied by Dixmier \cite{Dixmier.Masa} who used them to distinguish various types of maximal abelian subalgebras (masas).  He defined a masa $B\subseteq M$ to be singular, regular, or semiregular respectively if the normalising algebra is $B$, $M$, or a proper subfactor. For two inclusions $B_i\subseteq M_i$ of masas in finite von Neumann algebras, Chifan \cite{Chifan.Normalisers} showed that the normalising algebra of the tensor product is related to the tensor product of the individual normalising algebras by
\begin{equation}\label{NTensor}
\N_{M_1\,\vnotimes\,M_2}(B_1\ \vnotimes\ B_2)''=\N_{M_1}(B_1)''\ \vnotimes\ \N_{M_2}(B_2)'',
\end{equation}
a formula that was established earlier in \cite{Saw.StrongSing} for singular masas.  The commutation identity (\ref{NTensor}) also holds when each $B_i$ is an irreducible subfactor of a II$_1$ factor, \cite{Saw.NormalisersSubfactors}, and in this situation any normaliser $u$ is of the form $w(v_1\otimes v_2)$ for normalisers $v_i\in\N_{M_i}(B_i)$ and some unitary $w\in B_1\ \vnotimes\ B_2$. 

In general (\ref{NTensor}) fails. Consider the inclusion of $B=\mathbb C\oplus\mathbb M_2$ as a subalgebra of the $3\times 3$ matrices $\mathbb M_3=M$.   Dimension considerations show that every normaliser of $B$ lies in $B$. Since $B\otimes B\cong \mathbb C\oplus \mathbb M_2\oplus\mathbb M_2\oplus\mathbb M_4$, we can find a unitary in $\mathbb M_3\otimes\mathbb M_3$ which interchanges the two copies of $\mathbb M_2$.  This produces a non-trivial normaliser of the tensor product inclusion $B\otimes B\subseteq M\otimes M$.  The obstruction here is the presence of partial isometries $v\in M\setminus B$ with $vBv^*\subseteq B$ and $v^*Bv\subseteq B$.  Such partial isometries form a groupoid denoted $\GN_M(B)$ and are referred to as the \emph{groupoid normalisers} of $B$ in $M$. For masas in finite von Neumann algebras and irreducible inclusions of subfactors of II$_1$ factors the groupoid normalisers and normalisers generate the same von Neumann algebra and so this obstruction does not occur, but in general $\GN_M(B)''$ can be larger than $\N_M(B)''$. There are two related problems regarding the form of the groupoid normalising algebra of a tensor product of inclusions.
\begin{question}\label{MainQ}
Consider two unital inclusions $B_i\subseteq M_i$ of von Neumann algebras.
\begin{itemize}
\item[{\rm{(i)}}] Under what conditions on both inclusions do we have
\begin{equation}\label{GNCom}
\GN_{M_1\,\vnotimes\,M_2}(B_1\ \vnotimes\ B_2)''=\GN_{M_1}(B_1)''\ \vnotimes\ \GN_{M_2}(B_2)''?
\end{equation}
\item[{\rm{(ii)}}] Under what conditions on $B_1\subseteq M_1$ is (\ref{GNCom}) valid for all choices of $B_2\subseteq M_2$?
\end{itemize}
\end{question}

The main result of \cite{Saw.Groupoid} established (\ref{GNCom}) when $M_1$ and $M_2$ are finite and the relative commutants $B_i'\cap M_i$ lie in $B_i$ for $i=1,2$.    In general some assumption on the relative commutant is necessary for (\ref{GNCom}) to hold (see \cite[Example 1.1]{Saw.Groupoid} for a $3\times 3$ matrix example demonstrating this).  In the case when one of $M_1$ or $M_2$ is infinite, additional hypotheses will be needed to ensure that (\ref{GNCom}) holds.  Indeed, the following easy example shows that it is not possible to extend Chifan's result \cite{Chifan.Normalisers} for masas to the infinite setting. 
\begin{example}
Let $H=L^2[0,1]$ and let $A=L^\infty[0,1]$ acting on $H$ by left multiplication.  This is the unique (up to unitary conjugacy) diffuse masa in $\mathbb B(H)$ and direct computations show that $\N_{\mathbb B(H)}(A)''=\GN_{\mathbb B(H)}(A)''=\mathbb B(H)$.  Now let $H_1=H\oplus \mathbb C$ and let $A_1=A\oplus\mathbb C$. This is a masa in $\mathbb B(H_1)$ and 
\begin{equation}
\GN_{\mathbb B(H_1)}(A_1)''=\N_{\mathbb B(H_1)}(A_1)''=\mathbb B(H)\oplus \mathbb B(\mathbb C)\subsetneqq\mathbb B(H_1).
\end{equation}
Now consider $A\ \vnotimes\ A_1$ acting on $H\otimes H_1$. Since this is a diffuse masa, it is unitarily conjugate to the original inclusion $A\subseteq \mathbb B(H)$.  In particular
\begin{equation}
\GN_{\mathbb B(H)\,\vnotimes\,\mathbb B(H_1)}(A\ \vnotimes\ A_1)''=\N_{\mathbb B(H)\,\vnotimes\,\mathbb B(H_1)}(A\ \vnotimes\ A_1)''=\mathbb B(H)\ \vnotimes\ \mathbb B(H_1).
\end{equation}
Thus
\begin{equation}
\GN_{\mathbb B(H)\,\vnotimes\,\mathbb B(H_1)}(A\ \vnotimes\ A_1)''\supsetneqq\GN_{\mathbb B(H)}(A)''\ \vnotimes\ \GN_{\mathbb B(H_1)}(A_1)''.
\end{equation}
\end{example}

In this paper we obtain  positive answers to Question \ref{MainQ} (i) in several different contexts, and also to Question \ref{MainQ} (ii) in the following two situations:
\begin{itemize}
\item[{\rm{(1)}}] When one inclusion has an atomic relative commutant $B_1'\cap M_1$ which is equal to the centre of $B_1$ (Theorem \ref{Tensor.Main}).
\item[{\rm{(2)}}] When one inclusion is the generator masa in a free group factor (Theorem \ref{FreeGroup}).
\end{itemize}
Unlike \cite{Saw.Groupoid}, both cases above make no assumption on the other inclusion. The methods employed   there make extensive use of basic construction techniques and so are specialised to inclusions of finite von Neumann algebras.  Here we develop new techniques which reach outside the finite setting and so the paper can be read independently of \cite{Saw.Groupoid}. The second case described above is obtained by a direct combinatorial calculation examining the group elements supporting a groupoid normaliser.  This is a self-contained argument found in Lemma \ref{FreeGroupLem} and Theorem \ref{FreeGroup} below.  In particular, this gives the first  example of a  masa--factor inclusion $B_1\subseteq M_1$ with diffuse relative commutant $B_1'\cap M_1$ so that (\ref{GNCom}) holds for all inclusions $B_2\subseteq M_2$.

We study groupoid normalisers by examining the fixed point algebra of certain groups of automorphisms.  These techniques have their origins in \cite{Fang.MaxInjective,Fang.CompletelySingular} and are based on the simple observation that any normaliser $u\in\N_M(B)$ gives rise to an automorphism $\ad(u)$ of the commutant of $B$ (in some faithful representation) which fixes $M'$ pointwise.  In Section \ref{Aut} we examine all automorphisms of $B'$ which fix $M'$ pointwise and show that the subalgebra of $B'$ fixed by these automorphisms is precisely the commutant of the groupoid normalisers $\GN_M(B)$ (Theorem \ref{GN=N}).  This enables us to characterise those inclusions $B\subseteq M$ for which $\N_M(B)''=\GN_M(B)''$. In Section \ref{Tensor}, we use the description of $\GN_M(B)$ in terms of these automorphisms to establish instances of (\ref{GNCom}).    A key idea is to tensor by copies of $\mathbb B(H)$ to ensure that all our von Neumann algebras are properly infinite as, in the countably decomposable situation, this forces $\N_M(B)''=\GN_M(B)''$ (Lemma \ref{PropInfiniteGN=N}).

In Section \ref{Irred} we examine two inclusions of irreducible subfactors $B_i\subseteq M_i$.  In \cite{Saw.NormalisersSubfactors} it was shown that when both $M_1$ and $M_2$ are finite, then a normaliser of $B_1\ \vnotimes\ B_2$ in $M_1\ \vnotimes\ M_2$ must factorise as $w(v_1\otimes v_2)$ for $w\in B$ and $v_i\in\N_{M_i}(B_i)$.  Like \cite{Saw.Groupoid}, this relies on basic construction techniques.  In Section \ref{Irred} we study such normalisers $u$ by means of the induced automorphism $\ad(u)$ of the commutant $B_1'\ \vnotimes\ B_2'$.  We show that such an automorphism necessarily splits as a tensor product $\theta_1\otimes\theta_2$ of automorphisms of $B_i'$ which fix $M_i'$.  From this we recover \cite[Theorem 4.2]{Saw.NormalisersSubfactors} and extend it to cover the following situations:
\begin{itemize}
\item[{\rm{(1)}}] Only one of $M_1$ or $M_2$ is finite;
\item[{\rm{(2)}}] Both $B_1$ and $B_2$ are infinite;
\item[{\rm{(3)}}] $B_1$ is finite with trivial fundamental group and $M_1$ is infinite.
\end{itemize}
We provide examples that show that this result can fail whenever the algebra $B_1$ of (3) has a nontrivial fundamental group. This enables us to characterise II$_1$ factors with trivial fundamental group in terms of normalising unitaries of tensor product inclusions.

Many results in the paper rely on Lemma \ref{PropInfiniteTech} which requires a hypothesis of countable decomposability, and so attention is generally resticted to this class of algebras.
The exceptions are the results from Lemma \ref{FreeGroupLem} through to Theorem \ref{thm5.8}. Finally, all inclusions $B\subseteq M$ of von Neumann algebras in the paper are assumed to share the same unit  unless explicitly stated otherwise.  

\subsection*{Acknowledgements}
The work in this paper originated in the Workshop in Analysis and Probability, held at Texas A\&{}M in 2008. It is a pleasure to record our thanks to the organisers of the workshop and to the NSF for providing financial support to the workshop.

\section{Preliminaries}

In this section we establish some preliminary results, the most important of which is Lemma \ref{StdFix} which describes a useful connection between groupoid normalisers and the fixed point algebras of certain automorphism groups. We begin with the following observation, which is a version of \cite[Lemma 2.11]{Fang.CompletelySingular}. 
\begin{proposition}\label{Polar}
Let $B\subseteq M$ be an inclusion of von Neumann algebras and let $x\in M$ satisfy $xBx^*\subseteq B$ and $x^*Bx\subseteq B$.  Let $p,q\in B$ be  central projections and  let $pxq=v|pxq|$ be the polar decomposition of $pxq$.  Then $|pxq|\in B$ and $v\in \GN_M(B)$.
\end{proposition}

For a masa $A$ in a II$_1$ factor $M$, a theorem of Dye
 \cite{Dye.Groups2} (see also \cite[Lemma 6.2.3]{Sinclair.MasaBook}) states that  every $v\in\GN_M(A)$ is of the form $ue$ for some projection $e\in A$ and a unitary normaliser $u\in \N_M(A)$.  This cannot hold in the infinite setting: the unilateral shift is a groupoid normaliser of the   diagonal masa in $\mathbb B(\ell^2(\mathbb N))$ which cannot   be extended to a unitary normaliser.  However the unilateral shift does lie in the algebra generated by the normalisers. The lemma below expands this observation to   a general setting.  We adopt the notation $Z(N)$ for the centre of any von Neumann algebra $N$.

\begin{lemma}\label{CentreDye}
Let $B\subseteq M$ be an inclusion of von Neumann algebras and let $v\in\GN_M(B)$ satisfy $vv^*,v^*v\in Z(B)$. Then $v\in\N_M(B)''$.
\end{lemma}
\begin{proof}
Note that the hypotheses ensure that both $v$ and $v^*$ lie in $\GN_M(Z(B))$. Indeed, given $z\in Z(B)$ and $b\in B$, we have
\begin{equation}
vzv^*b=vzv^*bvv^*=vv^*bvzv^*=bvzv^*,
\end{equation}
so that $vZ(B)v^*\subseteq Z(B)$. Similarly $v^*Z(B)v^*\subseteq Z(B)$.

Let $p$ be the maximal projection in $Z(B)$ with $p\leq v^*v$ and $vp\in\N_M(B)''$. Suppose that $p\neq v^*v$ and write $p_0=v^*v-p\in Z(B)$ and $q_0=vp_0v^*\in Z(B)$.  If $p_0=q_0$, then $u=vp_0+(1-p_0)$ is a unitary normaliser of $B$ and $vp_0=up_0\in\N_M(B)''$, contradicting maximality of $p$. Hence $p_0\neq q_0$. Now suppose that $p_1=p_0(1-q_0)\neq 0$, and write $q_1=vp_1v^*\in Z(B)$ so $q_1p_0=0$.  Define a unitary by $u=vp_1+(vp_1)^*+(1-p_1-q_1)$. For $b\in B$, we have
\begin{equation}
ubu^*=u(p_1b+q_1b+(1-p_1-q_1)b)u^*=vp_1bv^*+v^*q_1bv+(1-p_1-q_1)b\in B
\end{equation}
so that $u$ is a unitary normaliser of $B$ with $vp_1=up_1$.  Thus $v(p+p_1)\in\N_M(B)''$, contradicting maximality of $p$.  Finally, if $q_1=(1-p_0)q_0\neq 0$, then by interchanging the roles of $v$ and $v^*$, there is some $u\in\N_M(B)$ with $u^*q_1=v^*q_1$ so that $q_1v=v(v^*(1-p_0)q_0v)\in\N_M(B)''$. As $v^*(1-p_0)q_0v\in Z(B)$,  we can adjoin $v^*(1-p_0)q_0v\leq p_0$ to $p$, contradicting the maximality of $p$.
\end{proof}

For inclusions of properly infinite von Neumann algebras, standard techniques allow us to adjust groupoid normalisers to have central initial and final projections so the previous lemma applies.  We record the details in slightly greater generality for use in Section \ref{Tensor}.
\begin{lemma}\label{PropInfiniteTech}
Let $M$ be a von Neumann algebra and let $B_1,B_2\subseteq M$ be properly infinite, countably decomposable von Neumann subalgebras of $M$.  Any partial isometry $v\in M$ with $vB_1v^*\subseteq B_2$ and $v^*B_2v\subseteq B_1$ factorises as $v=b_2wb_1$, where $b_1,b_2,w$ are partial isometries with $wB_1w^*\subseteq B_2$, $w^*B_2w\subseteq B_1$, $w^*w\in Z(B_1)$, $ww^*\in Z(B_2)$, $b_1\in B_1$ and $b_2\in B_2$.
\end{lemma}
\begin{proof}
Let $e=v^*v\in B_1$ and let $p$ be the central support of $e$ in $B_1$. As $B_1$ is properly infinite, standard arguments enable us to find a sequence $(e_n)_{n=1}^\infty$ (necessarily countable from the hypothesis of countable decomposability) of pairwise orthogonal, equivalent projections in $B_1$, and which sum to $p$ and with $e_1=e$. Similarly, let $vv^*=f\in B_2$, let $q$ be the central support of $f$ in $B_2$, and find a sequence $(f_n)_{n=1}^\infty$ of pairwise orthogonal, equivalent projections in $B_2$ which sum to $q$ with $f_1=f$. For each $n$, find partial isometries $b_{1,n}\in B_1$ with $b_{1,n}^*b_{1,n}=e$ and $b_{1,n}b_{1,n}^*=e_n$ and partial isometries $b_{2,n}\in B_2$ with $b_{2,n}^*b_{2,n}=f$ and $b_{2,n}b_{2,n}^*=f_n$.
Then define a partial isometry
\begin{equation}
w=\sum_{n=1}^\infty b_{2,n}vb_{1,n}^*\in M,
\end{equation}
noting that the series converges in the strong operator topology.  By construction $ww^*=q\in Z(B_2)$, $w^*w=p\in Z(B_1)$, $wB_1w^*\subseteq B_2$, $w^*B_2w\subseteq B_1$ and $v=b_{2,1}^*wb_{1,1}$.
\end{proof}

\begin{lemma}\label{PropInfiniteGN=N}
Let $B\subseteq M$ be an inclusion of von Neumann algebras. Suppose that $B$ is properly infinite and countably decomposable.  Then $\GN_M(B)''=\N_M(B)''$.
\end{lemma}
\begin{proof}
Using Lemma \ref{PropInfiniteTech}, any groupoid normaliser $v$ of $B$ can be factorised in the form $b_2wb_1$ so that $w$ is a groupoid normaliser of $B$ with $ww^*,w^*w\in Z(B)$ and $b_1,b_2\in B$. By Lemma \ref{CentreDye}, $w$ lies in $\N_M(B)''$ and hence so too does $v$. 
\end{proof}

Now we turn to the connections between normalisers, groupoid normalisers and certain automorphisms of the commutant inclusion.

\begin{definition}
Given an inclusion $B\subseteq M$ of von Neumann algebras,  we define the $B$-bimodular automorphisms of $M$, denoted $\Aut_B(M)$, to be those  $\theta\in\Aut(M)$ satisfying $\theta(b)=b$ for all $b\in B$.  Given a subgroup $G$ of the automorphism group of $M$, let $M^G$ denote the fixed point algebra $\{x\in M:\theta(x)=x,\ \theta\in G\}$.  This is a von Neumann subalgebra of $M$.  In particular we can apply this to $\Aut_B(M)$ to obtain $M^{\Aut_B(M)}$ which is a von Neumann subalgebra of $M$ containing $B$.

A number of the von Neumann algebras and groups that will appear below have complicated notations, so we will sometimes adopt the expression $\Fix(M,G)$ for $M^G$.
\end{definition}

On occasion we shall need to assume that certain von Neumann algebras are represented in \emph{standard form}.  In full generality, the standard form of a von Neumann algebra was set out by Haagerup, \cite{Haagerup.StandardForm}. The key fact we will need is that every automorphism of a von Neumann algebra in standard form is spatially implemented.
\begin{lemma}\label{StdFix}
Let $B\subseteq M \subseteq \mathbb{B}(H) $ be an inclusion of von Neumann algebras. Then 
\begin{equation}\label{StdFix.1}
(B')^{\Aut_{M'}(B')}\subseteq \N_M(B)',
\end{equation}
and equality holds if $B$ (or equivalently $B'$) is in standard form on $H$.
\end{lemma}

\begin{proof} 
Let $x\in(B')^{\Aut_{M'}(B')}$. Given $u\in\N_M(B)$, we have $uBu^*=B$ so $uB'u^*=B'$. Thus $\ad(u)|_{B'}$ defines an automorphism of $B'$ which is $M'$-bimodular as $u\in M$.  Hence $\ad(u)(x)=x$ and so $x\in \N_M(B)'$.

Now suppose that $B'$ lies in standard form on $H$ and take $x\in\N_M(B)'\subseteq B'$.  Given $\theta\in\Aut_{M'}(B')$,  there is a unitary $u\in\mathbb B(H)$ such that $\theta=\ad(u)|_{B'}$ as $B'$ is in standard form on $H$.   Since $\theta(m')=m'$ for $m'\in M'$,  the double commutant theorem gives $u\in M$. Now $uB'u^*=B'$, so we can take commutants to see that $uBu^*=B$, placing $u\in\N_M(B)$. It follows that $\theta(x)=uxu^*=x$, since $x\in N_M(B)'$, and so $x\in(B')^{\Aut_{M'}(B')}$, as required.
\end{proof}

\section{Groupoid normalisers and fixed point algebras}\label{Aut}

Our objective in this section is to characterise the groupoid normalising algebra $\GN_M(B)''$ of an inclusion of von Neumann algebras as the commutant of the fixed point algebra $(B')^{\Aut_{M'}(B')}$.  The first step is to show that the inclusions of an atomic masa inside $\mathbb B(K)$, and the inclusion $\mathbb B(K)\subseteq\mathbb B(K)$ give positive answers to Question \ref{MainQ} (ii). We proceed by direct computations with the matrix units as in \cite[Section 2.2]{Fang.CompletelySingular}.

\begin{lemma}\label{Atomic}
Let $K$ be a Hilbert space and let $A$ be the atomic masa in $\mathbb B(K)$. Given any inclusion $B\subseteq M$ of von Neumann algebras, we have
\begin{equation}\label{Atomic.1}
\GN_{M\,\vnotimes\,\mathbb B(K)}(B\ \vnotimes\ A)''=\GN_{M\,\vnotimes\,\mathbb B(K)}(B\ \vnotimes\ \mathbb B(K))''=\GN_M(B)''\ \vnotimes\ \mathbb B(K).
\end{equation}
\end{lemma}
\begin{proof}
We consider a more general situation which will imply containment of each of the first two algebras from (\ref{Atomic.1}) into the third simultaneously. Thus take
 a partial isometry $v\in M\ \vnotimes\ \mathbb B(K)$ with 
\begin{equation}\label{Atomic.2}
v(B\ \vnotimes\ A)v^*\subseteq B\ \vnotimes\ \mathbb B(K)\ {\mathrm{ and }}\ v^*(B\ \vnotimes\ A)v\subseteq B\ \vnotimes\ \mathbb B(K),
\end{equation}
and note that groupoid normalisers from the first two algebras in (\ref{Atomic.1}) will satisfy (\ref{Atomic.2}).
 Use the minimal projections in $A$ to identify operators in $M\ \vnotimes\ \mathbb B(K)$ with matrices over $M$, and write $v=(v_{i,j})_{i,j\in I}$.  For each $j\in I$ and $x\in B$, let $y\in B\ \vnotimes\ A$ be the operator with $x$ in the $(j,j)$-position and $0$ elsewhere.  By considering the $(i,i)$-th component of $vyv^*$, we see that $v_{i,j}xv_{i,j}^*\in B$.  Similarly $v_{i,j}^*Bv_{i,j}\subseteq B$ for all $i,j$.  Let $w_{i,j}h_{i,j}$ be the polar decomposition of $v_{i,j}$, so $w_{i,j}\in\GN_M(B)$ and $h_{i,j}\in B$ by Proposition \ref{Polar}.  Thus each $v_{i,j}$ lies in $\GN_M(B)''$ so that $v\in \GN_M(B)''\ \vnotimes\ \mathbb B(K)$.  This shows that the first and second algebras in (\ref{Atomic.1}) are contained in $\GN_M(B)''\ \vnotimes\ \mathbb B(K)$. Since the reverse inclusions are immediate, the result follows.
\end{proof}

\begin{lemma}\label{NotimesBH}
Let $B\subseteq M$ be an inclusion of von Neumann algebras such that $B$ is countably decomposable, and let $K$ be an infinite dimensional separable Hilbert space.  Then
\begin{equation}
\N_{M\,\vnotimes\,\mathbb B(K)}(B\ \vnotimes\ \mathbb B(K))''=\GN_M(B)''\ \vnotimes\ \mathbb B(K).
\end{equation}
\end{lemma}
\begin{proof}
Since $B\ \vnotimes\ \mathbb B(K)$ is properly infinite and countably decomposable, Lemma \ref{PropInfiniteGN=N} gives \begin{equation}
\N_{M\,\vnotimes\,\mathbb B(K)}(B\ \vnotimes\ \mathbb B(K))''=\GN_{M\,\vnotimes\,\mathbb B(K)}(B\ \vnotimes\ \mathbb B(K))''.
\end{equation}
 The result then follows from the second equality of Lemma \ref{Atomic}.
\end{proof}

We can now characterise those inclusions for which $\N_M(B)''=\GN_M(B)''$ in terms of the fixed points of the automorphism group $\Aut_{M'}(B')$.
\begin{theorem}\label{GN=N}
Let $B\subseteq M \subseteq \mathbb{B}(H) $ be an inclusion of von Neumann algebras such that $B$ is countably decomposable, and let $K$ be an infinite dimensional separable Hilbert space.  The following statements are equivalent:
\begin{enumerate}[(i)]
\item $\N_M(B)''=\GN_M(B)''$;\label{GN=N:L1}
\item $\N_{M\,\vnotimes\,\mathbb B(K)}(B\ \vnotimes\ \mathbb B(K))''=\N_M(B)''\ \vnotimes\ \mathbb B(K)$;\label{GN=N:L2}
\item $\N_M(B)'=(B')^{\Aut_{M'}(B')}$.\label{GN=N:L3}
\end{enumerate}
\end{theorem}
\begin{proof} The equivalence between statements (\ref{GN=N:L1}) and (\ref{GN=N:L2}) follows from Lemma \ref{NotimesBH} and is implicit in \cite[Theorem 2.4]{Fang.CompletelySingular}.

\medskip

\noindent $(\ref{GN=N:L3})\Longrightarrow (\ref{GN=N:L2})$.\quad Applying Lemma \ref{StdFix} to $B\ \vnotimes\ \mathbb B(K)\subseteq M\ \vnotimes\ \mathbb B(K)$ gives
\begin{equation}
\Fix (B'\otimes\mathbb C1,{\Aut_{M'\otimes \mathbb C1}(B'\otimes\mathbb C1)})\subseteq \N_{M\,\vnotimes\,\mathbb B(K)}(B\ \vnotimes\ \mathbb B(K))'.
\end{equation}
We have an isomorphism between $\Aut_{M'}(B')$ and $\Aut_{M'\otimes \mathbb C1}(B'\otimes\mathbb C1)$ given by  $\theta\mapsto\theta\otimes \mathrm{id}$  which demonstrates that
\begin{equation}
\Fix (B'\otimes\mathbb C1,{\Aut_{M'\otimes \mathbb C1}(B'\otimes\mathbb C1)})=(B')^{\Aut_{M'}(B')}\otimes\ \mathbb C1.
\end{equation}
Applying the hypothesis of (3) gives
\begin{equation}\label{eq3.5}
\N_M(B)'\otimes\mathbb C1\subseteq\N_{M\,\vnotimes\,\mathbb B(K)}(B\ \vnotimes\ \mathbb B(K))',
\end{equation}
and the inclusion
\begin{equation}
\N_{M\,\vnotimes\,\mathbb B(K)}(B\ \vnotimes\ \mathbb B(K))''\subseteq\N_M(B)''\ \vnotimes\ \mathbb B(K)''
\end{equation}
follows by taking commutants in (\ref{eq3.5}).
Since the reverse inclusion is immediate, condition (2) holds.

\medskip

\noindent $(2)\Longrightarrow (3)$.\quad  Let $\pi$ be a standard representation of $B'$ on some Hilbert space $H_1$.  By the general theory of representations of von Neumann algebras \cite[p. 61]{Dixmier.vNaBook}, we can assume that $\pi$ is obtained by an amplification of the representation on $H$ followed by a compression.  Therefore, we can find another Hilbert space $K$ and a projection in $e\in B\ \vnotimes\ \mathbb B(K)=(B'\otimes 1_K)'$ such that $\pi(x)=e(x\otimes 1_K)e$ acting on $e(H\otimes K)$.  Writing $B_1=\pi(B')'$ and $M_1=\pi(M')'$ with the commutants taken on $H_1=e(H\otimes K)$, we have $B_1=e(B\ \vnotimes\ \mathbb B(K))e$ and $M_1=e(M\  \vnotimes\ \mathbb B(K))e$.  Since $\pi$ is a faithful representation of $B'$, we have $\pi((B')^{\Aut_{M'}(B')})=(B_1')^{\Aut_{M_1'}(B_1')}$. Furthermore, since $\pi(B')$ is in standard form, so too is $B_1=\pi(B')'$. Thus Lemma \ref{StdFix} gives 
\begin{equation}\label{GN=N.3}
\pi((B')^{\Aut_{M'}(B')})=(B_1')^{\Aut_{M_1'}(B_1')}=\N_{M_1}(B_1)'.
\end{equation}

Every normaliser $v$ of $e(B\ \vnotimes\ \mathbb B(K))e$ in $e(M\ \vnotimes\ \mathbb B(K))e$ is a groupoid normaliser of $B\ \vnotimes\ \mathbb B(K)$ in $M\ \vnotimes\ \mathbb B(K)$ and so lies in $\N_{M\,\vnotimes\,\mathbb B(K)}(B\ \vnotimes\ \mathbb B(K))''$ by Lemma \ref{PropInfiniteGN=N}.  This gives the inclusion
\begin{equation}\label{GN=N.2}
\N_{M_1}(B_1)''=\N_{e(B\,\vnotimes\,\mathbb B(K))e}(e(M\ \vnotimes\ \mathbb B(K))e)''\subseteq e(\N_{M\,\vnotimes\,\mathbb B(K)}(B\ \vnotimes\ \mathbb B(K))'')e.
\end{equation}
Now 
\begin{equation}\label{GN=N.1}
\pi(\N_M(B)')'=[(\N_M(B)'\ \vnotimes\ \mathbb C1_{K})e]'=e(\N_M(B)''\ \vnotimes\ \mathbb B(K))e,
\end{equation}
so that 
\begin{align}
\pi(\N_M(B)')'&=e(\N_M(B)''\ \vnotimes\ \mathbb B(K))e\notag\\
&=e(\N_{M\,\vnotimes\,\mathbb B(K)}(B\ \vnotimes\ \mathbb B(K))'')e\notag\\
&\supseteq \N_{e(B\,\vnotimes\,\mathbb B(K))e}(e(M\ \vnotimes\ \mathbb B(K))e)''\notag\\
&=\N_{M_1}(B_1)''=\pi((B')^{\Aut_{M'}(B')})',
\end{align}
by (\ref{GN=N.1}), condition (2), (\ref{GN=N.2}) and (\ref{GN=N.3}). Taking commutants gives 
\begin{equation}
\pi(\N_M(B)')\subseteq\pi((B')^{\Aut_{M'}(B')}),
\end{equation}
and so 
\begin{equation}
\N_M(B)'\subseteq(B')^{\Aut_{M'}(B')}
\end{equation}
since $\pi$ is a faithful representation of $B'$.  The reverse inclusion $(B')^{\Aut_{M'}(B')}\subseteq\N_M(B)'$ is Lemma \ref{StdFix} so condition (3) holds.
\end{proof}

Finally in this section we can express the groupoid normalising algebra $\GN_M(B)''$ in terms of the automorphism group $\Aut_{M'}(B')$.  We will use this result repeatedly in the next section to obtain further instances of (\ref{GNCom}).
\begin{theorem}\label{GN=Fix}
Let $B\subseteq M\subseteq \mathbb B(H)$ be an inclusion of von Neumann algebras where $B$ is countably decomposable.  Then 
\begin{equation}
\GN_M(B)'=(B')^{\Aut_{M'}(B')}.
\end{equation}
\end{theorem}
\begin{proof}
Let $K$ be an infinite dimensional separable Hilbert space and define $B_1=B\ \vnotimes\ \mathbb B(K)$ and $M_1=M\ \vnotimes\ \mathbb B(K)$.  Then
\begin{equation}\label{GN=Fix.1}
\N_{M_1}(B_1)''=\GN_{M}(B)''\ \vnotimes\ \mathbb B(K), 
\end{equation}
by Lemma \ref{NotimesBH}. Since $B_1$ is properly infinite, Lemma \ref{PropInfiniteGN=N} gives $\GN_{M_1}(B_1)''=\N_{M_1}(B_1)''$ so the inclusion $B_1\subseteq M_1$ satisfies condition (1) of Theorem \ref{GN=N}. By (\ref{GN=Fix.1}) and condition (3) of this theorem,
\begin{equation}
\GN_M(B)' \otimes \mathbb C1_{K}=\N_{M_1}(B_1)'=(B_1')^{\Aut_{M_1'}(B_1')}.
\end{equation}
Finally note that $\theta\mapsto \theta\,\otimes\,\mathrm{id}_{\mathbb B(K)}$ gives an isomorphism between $\Aut_{M'}(B')$ and $\Aut_{M_1'}(B_1')$ so that $(B_1')^{\Aut_{M_1'}(B_1')}=(B')^{\Aut_{M'}(B')}\ \vnotimes\ \mathbb C1_K$. Thus $\GN_M(B)'=(B')^{\Aut_{M'}(B')}$.
\end{proof}

\section{Groupoid normalisers and tensor products}\label{Tensor}

In this section, we establish our positive answers to Question \ref{MainQ}. We being with separably acting inclusions of the form $B\ \vnotimes\ A\subseteq M\ \vnotimes\ A$, where $A$ is abelian. When $B$ is singular in $M$, \cite[Lemma 4.1]{Fang.CompletelySingular} showed that $B\ \vnotimes\ A$ is singular in $M\ \vnotimes\ A$ by a direct integral argument based on \cite[Lemma 6.6]{Stratila.Commutation}. The proof given in \cite{Fang.CompletelySingular} shows that any normaliser $u$ of $B\ \vnotimes\ A$ in $M\ \vnotimes\ A$ is a direct integral of normalisers of $B$ in $M$ over the spectrum of $A$.  The analogous result also holds for groupoid normalisers, giving the following lemma. We omit the proof, which is a routine modification of \cite[Lemma 6.6]{Stratila.Commutation} and \cite[Lemma 4.1]{Fang.CompletelySingular}.
\begin{lemma}\label{AbDirectInt}
Let $B\subseteq M$ be an inclusion of separably acting von Neumann algebras.  If $A$ is a separable abelian von Neumann algebra, then 
\begin{equation}
\N_{M\,\vnotimes\,A}(B\ \vnotimes\ A)''=\N_M(B)''\ \vnotimes\ A,\quad\GN_{M\,\vnotimes\,A}(B\ \vnotimes\ A)''=\GN_M(B)''\ \vnotimes\ A.
\end{equation}
\end{lemma}

The next lemma follows the same pattern as the deduction of 
\cite[Theorem 4.3]{Fang.CompletelySingular} from \cite[Theorem 3.1, Lemma 4.1]{Fang.CompletelySingular}. 

\begin{lemma}\label{otimesL}
Let $B\subseteq M$ be a unital inclusion of separably acting von Neumann algebras and let $L$ be a separably acting von Neumann algebra.  Then
\begin{equation}
\GN_{M\,\vnotimes\,L}(B\ \vnotimes\ L)''=\GN_M(B)''\ \vnotimes\ L.
\end{equation}
\end{lemma}
\begin{proof}
Fix a faithful representation of $M$ on a separable Hilbert space $H_1$ and use this to define the commutants $B'\supseteq M'$. Write $A=Z(L)$ and take a standard representation of $A$ on a separable Hilbert space $K$ so that $A$ is a masa in $\mathbb B(K)$.  Work on the Hilbert space $H_1\otimes K$ so that $(B\ \vnotimes\ A)'=B'\ \vnotimes\ A$.   Given $\theta\in\Aut_{M'\,\vnotimes\,A}(B'\ \vnotimes\ A)$, Theorem \ref{GN=Fix} gives $\theta(x)=x$ for all $x\in\GN_{M\,\vnotimes\,A}(B\ \vnotimes\ A)'=\GN_{M}(B)'\ \vnotimes\ A$, where the last identity is obtained by taking commutants in Lemma \ref{AbDirectInt}.

Now suppose  that $L$ is faithfully represented on a separable Hilbert space $H_2$ and work on $H_1\otimes H_2$.  Given an automorphism $\theta\in\Aut_{M'\,\vnotimes\,L'}(B'\ \vnotimes\ L')$,   we have
\begin{align}
\theta(x)(1_{B'}\otimes \ell')&=\theta(x)\theta(1_{B'}\otimes \ell')=\theta(x(1_{B'}\otimes \ell'))=\theta((1_{B'}\otimes \ell')x)\notag\\
&=\theta(1_{B'}\otimes \ell')\theta(x)=(1_{B'}\otimes \ell')\theta(x),\ \ \ell'\in L',\ \  x\in B'\,\vnotimes \, A, 
\end{align}
since $\theta(1_{B'}\otimes \ell')=1_{B'}\otimes \ell'$ and $1_{B'}\otimes \ell'$ commutes with $B'\ \vnotimes\ A\subseteq B'\ \vnotimes\ L'$.  Thus 
\begin{equation}
\theta(x)\in (\mathbb C1_{B'}\otimes L')'\cap (B'\ \vnotimes\ L')=B'\ \vnotimes\ Z(L)=B'\ \vnotimes\ A,
\end{equation}
and so every element $\theta$ of $\Aut_{M'\,\vnotimes\,L'}(B'\ \vnotimes\ L')$ restricts to an element of $\Aut_{M'\,\vnotimes\,A}(B'\ \vnotimes\ A)$.  It then follows from the first paragraph of the proof that $\theta(x)=x$ for all $x\in\GN_M(B)'\ \vnotimes\ A$. Since $\theta$ was arbitrary, it follows that
\begin{equation}
\GN_M(B)'\ \vnotimes\ A\subseteq\Fix(B'\ \vnotimes\ L',\Aut_{M'\,\vnotimes\,L'}(B'\ \vnotimes\ L'))=\GN_{M\,\vnotimes\,L}(B\ \vnotimes\ L)',
\end{equation}
where the second equality is Theorem \ref{GN=Fix}. Take commutants to obtain
\begin{equation}
\GN_M(B)''\ \vnotimes\ A'\supseteq \GN_{M\,\vnotimes\,L}(B\ \vnotimes\ L)''
\end{equation}
on $H_1\otimes H_2$. Then
\begin{equation}
\GN_{M\,\vnotimes L}(B\ \vnotimes\ L)''\subseteq (\GN_M(B)''\ \vnotimes\ A')\cap (M\ \vnotimes\ L)=\GN_M(B)''\ \vnotimes\ L.
\end{equation}
Since the reverse inclusion is immediate, the result follows.
\end{proof}

We now start work on our main result. Given two unital inclusions $B_i\subseteq M_i$ of von Neumann algebras, we will show that 
\begin{equation}\label{eq4.8}
\GN_{M_1\,\vnotimes\,M_2}(B_1\ \vnotimes\ B_2)''=\GN_{M_1}(B_1)''\ \vnotimes\ \GN_{M_2}(B_2)'',
\end{equation}
when $B_1'\cap M_1$ is atomic and lies in the centre of $B_1$.   The inclusion from right to left is immediate. We establish the inclusion from left to right by demonstrating that $\GN_{M_1\,\vnotimes\,M_2}(B_1\ \vnotimes\ B_2)''$ is contained in $M_1\ \vnotimes\ \GN_{M_2}(B_2)''$ and in $\GN_{M_1}(B_1)''\ \vnotimes\ M_2$ separately.  The next lemma, which is based on \cite[Theorem 4.4]{Fang.CompletelySingular}, handles the first, and easier, of these two inclusions. \begin{lemma}\label{4.3}
For $i=1,2$, let $B_i\subseteq M_i$ be unital inclusions of separably acting von Neumann algebras and suppose that $B_1'\cap M_1=Z(B_1)$ is atomic.  Then
\begin{equation}
\GN_{M_1\,\vnotimes\,M_2}(B_1\ \vnotimes\ B_2)''\subseteq M_1\ \vnotimes\ \GN_{M_2}(B_2)''.
\end{equation}
\end{lemma}

\begin{proof}
Represent $M_1$ and $M_2$ faithfully on Hilbert spaces $H_1$ and $H_2$ respectively and consider an automorphism $\theta\in \Aut_{M_1'\,\vnotimes\,M_2'}(B_1'\ \vnotimes \ B_2')$.  Write $A=Z(B_1)$. For $y_1\in A=Z(B_1')$ and $y_2\in B_2'$, we have
\begin{align}
\theta(y_1\otimes y_2)(x\otimes 1)&=\theta(y_1x\otimes y_2)=\theta(xy_1\otimes y_2)\notag\\
&=(x\otimes 1)\theta(y_1\otimes y_2),\quad x\in M_1'.
\end{align}
Thus 
\begin{equation}
\theta(y_1\otimes y_2)\in (M_1'\otimes \mathbb C1)'\cap (B_1'\ \vnotimes\ B_2')=(B_1'\cap M_1)\ \vnotimes\ B_2'=A\ \vnotimes\ B_2',
\end{equation}
and so it follows that $\theta$ restricts to an element of $\Aut_{\mathbb C1\otimes M_2'}(A\ \vnotimes\ B_2')$.  Let $A$ be represented as the diagonal operators on some Hilbert space $K$, so that working on $K\otimes H_2$ we have $(A\ \vnotimes\ B_2)'=A\ \vnotimes\ B_2'$. Lemma \ref{Atomic} gives $\GN_{\mathbb B(K)\,\vnotimes\,M_2}(A\ \vnotimes\ B_2)''=\mathbb B(K)\ \vnotimes\ \GN_{M_2}(B_2)''$ so that 
\begin{equation}
\Fix(A\ \vnotimes\ B_2',\Aut_{\mathbb C1\otimes M_2'}(A\ \vnotimes\ B_2'))=\mathbb C1\otimes \GN_{M_2}(B_2)',
\end{equation}
by Theorem \ref{GN=Fix}. Thus $\theta(1\otimes z)=1\otimes z$ for all $z\in \GN_{M_2}(B_2)'$.  Theorem \ref{GN=Fix} also gives
\begin{equation}
\GN_{M_1\,\vnotimes\,M_2}(B_1\ \vnotimes\ B_2)'=\Fix(B_1'\ \vnotimes\ B_2',\Aut_{M_1'\,\vnotimes\,M_2'}(B_1'\ \vnotimes\ B_2'))
\end{equation}
so that 
\begin{equation}
\mathbb C1\otimes\GN_{M_2}(B_2)'\subseteq\GN_{M_1\,\vnotimes\,M_2}(B_1\ \vnotimes\ B_2)'.
\end{equation}
The inclusions
\begin{align}
\GN_{M_1\,\vnotimes\,M_2}(B_1\ \vnotimes\ B_2)''&\subseteq \left(\mathbb B(H_1)\ \vnotimes\ \GN_{M_2}(B_2)''\right)\cap(M_1\ \vnotimes\ M_2)\notag\\
&=M_1\ \vnotimes\ \GN_{M_2}(B_2)''
\end{align}
follow by taking commutants.
\end{proof}

Now we turn to the second inclusion required for (\ref{eq4.8}). Recall that a von Neumann algebra in standard form has the property that its automorphisms are all spatially implemented.

\begin{lemma}\label{stdformlem}
Let $B\subseteq M$ be a unital inclusion of separably acting von Neumann algebras with $B'\cap M\subseteq B$ and suppose that $p$ and $q$ are minimal central projections in $B$.  Let $L$ be another separable von Neumann algebra  acting in standard form on the Hilbert space $K$.  Suppose that $v\in\GN_{M\,\vnotimes\,\mathbb B(K)}(B\ \vnotimes\ L)$ satisfies $v^*v,vv^*\in Z(B\ \vnotimes\ L)$ with $v^*v\leq p\otimes 1$ and $vv^*\leq q\otimes 1$.  Then $v\in \GN_M(B)''\ \vnotimes \ \mathbb B(K)$.
\end{lemma}

\begin{proof}
We represent $M$ on a Hilbert space $H$.  Since $p$ and $q$ are minimal central projections in $Z(B)$, there are central projections $p_1,q_1\in L$ with $v^*v=p\otimes p_1$ and $vv^*=q\otimes q_1$. Thus $x\mapsto vxv^*$ gives an isomorphism from $(Bp)\ \vnotimes\ (Lp_1)$ onto $(Bq)\ \vnotimes\ (Lq_1)$.  As a consequence, we obtain a surjective isomorphism $\theta:(B'p)\ \vnotimes\ (L'p_1)\rightarrow (B'q)\ \vnotimes\ (L'q_1)$ by $\theta(y)=vyv^*$.  Note that $v$ commutes with $M'\otimes\mathbb C1$ so 
\begin{equation}\label{4.4.1}
\theta(M'\otimes p_1)=v(M'\otimes p_1)v^*=v(M'\otimes 1)v^*=(M'\otimes 1)vv^*=M'q\otimes q_1.
\end{equation}
We claim that $\theta(p\otimes L'p_1)=q\otimes L'q_1$.  Indeed, for $m'\in M'$ and $\ell'\in L'$, 
\begin{equation}
(m'p\otimes p_1)(p\otimes \ell'p_1)=(p\otimes \ell'p_1)(m'p\otimes p_1),
\end{equation}
so that applying $\theta$ and using (\ref{4.4.1}), we obtain
\begin{equation}
\theta(p\otimes \ell'p_1)\in (M'q\otimes q_1)'\cap (B'q\otimes L'q_1)=(B'\cap M)q\otimes L'q_1=q\otimes L'q_1,
\end{equation}
from the minimality of $q\in Z(B)$. Interchanging the roles of $v$ and $v^*$ shows that $\theta^{-1}(q\otimes L'q_1)=p\otimes L'p_1$, establishing the claim. 

Since both $L'p_1$ and $L'q_1$ act in standard form on $p_1(K)$ and $q_1(K)$ respectively, there is a partial isometry $w\in\mathbb B(K)$ with $w^*w=p_1$ and $ww^*=q_1$ so that 
\begin{equation}
v(p\otimes \ell'p_1)v^*=\theta(p\otimes \ell'p_1)=(1\otimes w)(q\otimes\ell'p_1)(1\otimes w)^*,\quad \ell'\in L'.
\end{equation}
Define $v_1=(1\otimes w)^*v$ so that
\begin{equation}
v_1v_1^*=(1\otimes w)^*vv^*(1\otimes w)=(1\otimes w)^*(q\otimes q_1)(1\otimes w)=q\otimes p_1
\end{equation}
and
\begin{equation}
v_1^*v_1=v^*(1\otimes w)(1\otimes w)^*v=v^*(1\otimes q_1)v=p\otimes p_1.
\end{equation}
As $(1\otimes w)(1\otimes w^*)=1\otimes q_1\geq v_1v_1^*$, we have $v=(1\otimes w)v_1$.  Certainly $v_1\in M\otimes p_1\mathbb B(K)p_1$, and direct computations give
\begin{align}
v_1(1\otimes \ell')=&v_1v_1^*v_1(1\otimes \ell')v_1^*v_1=v_1(p\otimes l'p_1)v_1^*v_1\nonumber\\
=&(1\otimes w)^*\theta(p\otimes \ell'p_1)(1\otimes w^*)v_1\nonumber\\
=&(q\otimes w^*w\ell'w^*w)v_1=(1\otimes \ell')v_1,\quad \ell'\in L'.\label{4.4.2}
\end{align}
Thus 
\begin{equation}
v_1\in (\mathbb C1\otimes L')'\cap (M\ \vnotimes\ p_1\mathbb B(K)p_1)=M\ \vnotimes\ Lp_1\subseteq M\ \vnotimes\ L.
\end{equation}
Since $x\mapsto v_1xv_1^*$ is an isomorphism from $B'p\ \vnotimes\ L'p_1 =(B\ \vnotimes\ L)'\cap \mathbb B(pH\otimes p_1K)$ onto $B'q\ \vnotimes\ L'p_1=(B\ \vnotimes\ L)'\cap \mathbb B(qH\otimes p_1K)$ we can take commutants to see that
\begin{equation}
v_1(B\ \vnotimes\ L)v_1^*\subseteq B\ \vnotimes\ L
\end{equation}
and similarly that $v_1^*(B\ \vnotimes\ L)v_1^*\subseteq B\ \vnotimes\ L$.  Consequently $v_1\in \GN_{M\,\vnotimes\,L}(B\ \vnotimes\ L)$ and so $v_1\in\GN_{M}(B)''\ \vnotimes\ L$ by Lemma \ref{otimesL}. It follows that
\begin{equation}
v=(1\otimes w)v_1\in \GN_{M}(B)''\ \vnotimes\ \mathbb B(K),
\end{equation}
as required.
\end{proof}

\begin{lemma}\label{4.5}
Let $B\subseteq M$ be an inclusion of separably acting von Neumann algebras so that $B'\cap M=Z(B)$ is atomic.  Let $L$ be another von Neumann algebra, acting in standard form on a separable Hilbert space $K$.  Then
\begin{equation}\label{4.5.1}
\GN_{M\,\vnotimes\,\mathbb B(K)}(B\ \vnotimes\ L)''\subseteq\GN_M(B)''\ \vnotimes\ \mathbb B(K).
\end{equation}
\end{lemma}
\begin{proof}
We first show that it suffices to prove the lemma under the additional assumption that $B$ is properly infinite.  If $B$ is not properly infinite, consider the inclusion 
\begin{equation}
B_0=\mathbb B(\ell^2(\mathbb N))\ \vnotimes\ B\subseteq \mathbb B(\ell^2(\mathbb N))\ \vnotimes\ M=M_0,
\end{equation}
 which also satisfies the hypotheses of the lemma and $B_0$ is properly infinite.  Lemma \ref{Atomic} gives
\begin{equation}
\GN_{M_0}(B_0)''= \mathbb B(\ell^2(\mathbb N))\ \vnotimes\ \GN_M(B)''
\end{equation}
and
\begin{equation}
 \GN_{M_0\,\vnotimes\,\mathbb B(K)}(B_0\ \vnotimes\ L)''=\mathbb B(\ell^2(\mathbb N))\ \vnotimes\ \GN_{M\,\vnotimes\,\mathbb B(K)}(B\ \vnotimes\ L)''.
\end{equation}
Thus, assuming that the lemma holds for the inclusion $B_0\subseteq M_0$, we obtain
\begin{align}
\mathbb B(\ell^2(\mathbb N))\ \vnotimes\ \GN_{M\,\vnotimes\,\mathbb B(K)}(B\ \vnotimes\ L)''=&\,\GN_{M_0\,\vnotimes\,\mathbb B(K)}(B_0\ \vnotimes\ L)''\notag\\
\subseteq&\,\GN_{M_0}(B_0)''\ \vnotimes\ \mathbb B(K)\notag\\
=&\,\mathbb B(\ell^2(\mathbb N))\ \vnotimes\ \GN_M(B)''\ \vnotimes\ \mathbb B(K),
\end{align}
and so (\ref{4.5.1}) holds for the original inclusion $B\subseteq M$.

Now assume that $B$ is properly infinite. Fix $v\in\GN_{M\,\vnotimes\,\mathbb B(K)}(B\ \vnotimes\ L)$ and minimal projections $p,q\in Z(B)$. We will show that $(p\otimes 1)v(q\otimes 1)\in\GN_M(B)''\ \vnotimes\ \mathbb B(K)$, from which the result  follows immediately. Let $(p\otimes 1)v(q\otimes 1)=w|(p\otimes 1)v(1\otimes 1)q|$ be the polar decomposition of $(p\otimes 1)v(q\otimes 1)$, so that $|(p\otimes 1)v(q\otimes 1)|\in B\ \vnotimes\ L$ and $w\in\GN_{M\,\vnotimes\,\mathbb B(K)}(B\ \vnotimes\ L)$ by Proposition \ref{Polar}.  We must show that $w\in\GN_{M}(B)''\ \vnotimes\ \mathbb B(K)$. Write $A_p=(p\otimes 1)(B\ \vnotimes\ L)(p\otimes 1)$ and $A_q=(q\otimes 1)(B\ \vnotimes\ L)(q\otimes 1)$ so that $wA_qw^*\subseteq A_p$ and $w^*A_pw\subseteq A_q$.  Note that both $A_p$ and $A_q$ are central cutdowns of $B\ \vnotimes\ L$ and so properly infinite. By Lemma \ref{PropInfiniteTech}, we can factorise $w=a_qua_p$, where $a_p$ and $a_q$ are partial isometries in $A_p$ and $A_q$ respectively and $u$ is  a partial isometry with $uA_qu^*\subseteq A_p$, $u^*A_pu\subseteq A_q$ and $u^*u\in Z(A_q)$, $uu^*\in Z(A_p)$.  Lemma \ref{stdformlem} then shows that $u$, and hence $w$, lies in $\GN_{M}(B)''\ \vnotimes\ \mathbb B(K)$, exactly as required.
\end{proof}

We can now establish (\ref{eq4.8}), giving a general class of inclusions with a positive answer to Question \ref{MainQ} (ii).
\begin{theorem}\label{Tensor.Main}
Let $B_1\subseteq M_1$ be an inclusion of separably acting von Neumann algebras where $B_1'\cap M_1=Z(B_1)$ is atomic and let $B_2\subseteq M_2$ be another unital inclusion of 
separably acting von Neumann algebras.  Then
\begin{equation}
\GN_{M_1\,\vnotimes\,M_2}(B_1\ \vnotimes\ B_2)''=\GN_{M_1}(B_1)''\ \vnotimes\ \GN_{M_2}(B_2)''.
\end{equation}
\end{theorem}
\begin{proof}
Take faithful representations of $M_1$ and $M_2$ on Hilbert spaes $H_1$ and $H_2$ respectively. Fix $\theta\in\Aut_{M_1'\,\vnotimes\,M_2'}(B_1'\ \vnotimes\ B_2')$. We will show that $\theta(x\otimes 1)=x\otimes 1$ for all $x\in\GN_{M_1}(B_1)'$.

 Define $L=B_2'\cap M_2$.  For $m_2'\in M_2'$, $b_1'\in B_1'$ and $\ell\in L=B_2'\cap M_2$ we have
\begin{equation}
\theta(b_1'\otimes \ell)(1\otimes m_2')=\theta(b_1'\otimes \ell m_2)=\theta(b_1'\otimes m_2\ell )=(1\otimes m_2')\theta(b_1'\otimes \ell ),
\end{equation}
so that $\theta(b_1'\otimes \ell )\in (\mathbb C1\otimes M_2')'\cap (B_1'\ \vnotimes\ B_2')=B_1'\ \vnotimes\ L$.  In particular $\theta$ restricts to an element of $\Aut_{M_1'\otimes \mathbb C1}(B_1'\ \vnotimes\ L)$.  Now take a standard representation of $L$ on $K$ and work on $H_1\otimes K$.  Applying Theorem \ref{GN=Fix} to the inclusion $B_1\ \vnotimes\ L'\subseteq M_1\ \vnotimes\ \mathbb B(K)$, we have
\begin{equation}
\Fix(B_1'\ \vnotimes\ L,\Aut_{M_1'\otimes \mathbb C1}(B_1'\ \vnotimes\ L))=\GN_{M_1\,\vnotimes\,\mathbb B(K)}(B\ \vnotimes\ L')'.
\end{equation}
In particular $\theta(z)=z$ for all $z\in\GN_{M_1\,\vnotimes\,\mathbb B(K)}(B\ \vnotimes\ L')'$.  Lemma \ref{4.5} gives
\begin{equation}\label{eq4.1}
\GN_{M_1\,\vnotimes\,\mathbb B(K)}(B_1\ \vnotimes\ L')''\subseteq\GN_{M_1}(B_1)''\ \vnotimes\ \mathbb B(K),
\end{equation}
and the inclusion
\begin{equation}
\GN_{M_1}(B_1)'\otimes\mathbb C1\subseteq\GN_{M_1\,\vnotimes\,\mathbb B(K)}(B_1\ \vnotimes\ L')'
\end{equation}
follows by taking commutants in (\ref{eq4.1}).
Thus $\theta(x\otimes 1)=x\otimes 1$ for all $x\in\GN_{M_1}(B_1)'$.

Applying Theorem \ref{GN=Fix} to the original inclusion $B_1\ \vnotimes\ B_2\subseteq M_1\ \vnotimes\ M_2$ gives
\begin{equation}
\GN_{M_1\,\vnotimes\,M_2}(B_1\ \vnotimes\ B_2)'=\Fix(B_1'\ \vnotimes\ B_2',\Aut_{M_1'\,\vnotimes\,M_2'}(B_1'\ \vnotimes\ B_2'))
\end{equation}
and so 
\begin{equation}
\GN_{M_1}(B_1)'\otimes \mathbb C1\subseteq \GN_{M_1\,\vnotimes\,M_2}(B_1\ \vnotimes\ B_2)'.
\end{equation}
Taking commutants gives
\begin{equation}
\GN_{M_1\,\vnotimes\,M_2}(B_1\ \vnotimes\ B_2)''\subseteq\GN_{M_1}(B_1)''\ \vnotimes\ \mathbb B(H_2).
\end{equation}
As Lemma \ref{4.3} states that
\begin{equation}
\GN_{M_1\,\vnotimes\,M_2}(B_1\ \vnotimes\ B_2)''\subseteq M_1\ \vnotimes\ \GN_{M_2}(B_2)'',
\end{equation}
we obtain
\begin{align}
\GN_{M_1\,\vnotimes\,M_2}(B_1\ \vnotimes\ B_2)''&\subseteq \left(\GN_{M_1}(B_1)''\ \vnotimes\ \mathbb B(H_2)\right)\cap \left(M_1\ \vnotimes\ \GN_{M_2}(B_2)''\right)\nonumber\\
&=\GN_{M_1}(B_1)''\ \vnotimes\ \GN_{M_2}(B_2)''.
\end{align}
As the reverse inclusion is immediate, the result follows.
\end{proof}

We end this section by giving an example of a diffuse masa-factor inclusion $A\subseteq M$ for which 
\begin{equation}\label{Gen.1}
\GN_{M\,\vnotimes\,N}(A\ \vnotimes\ B)''=\GN_{M}(A)''\ \vnotimes\ \GN_{N}(B)'',
\end{equation}
for all inclusions $B\subseteq N$.  Fix $k\in\mathbb N$ with $k\geq 2$ and let $\mathbb F_k$ be the free group on the $k$ generators $a,b_1,\dots,b_{k-1}$ (the argument given below will also be valid for $k=\infty$).  Let $M=L(\mathbb F_k)$ be the group von   Neumann algebra generated by $\mathbb F_k$.  We identify elements of the group with the corresponding elements in $M$ and let $A$ be the von Neumann subalgebra of $M$ generated by $a$.  This is a masa in $M$, known as a generator masa so has $A'\cap M=A$.  The subgroup generated by $a$ satisfies Dixmier's criterion for  singularity of the  masa $A$ \cite{Dixmier.Masa}   (see also \cite[p. 22]{Sinclair.MasaBook}) and so  $\GN_M(A)''=A$. The inclusion $A\subseteq M$ fits into the framework of \cite{Saw.Groupoid} so (\ref{Gen.1}) holds whenever $B\subseteq N$ is an inclusion of finite von Neumann algebras with $B'\cap N\subseteq B$.  Our objective is to establish (\ref{Gen.1}) without making any assumption on the inclusion $B\subseteq N$.

We denote the standard  orthonormal basis for $\ell^2(\mathbb F_k)$ by $\{\delta_g:g\in \mathbb F_k\}$. When we view an operator $x\in L(\mathbb F_k)$ as a vector in the underlying Hilbert space, then it has a square summable Fourier series $\sum_{g\in \mathbb F_k}\alpha_g\delta_g$. The support of $x$ is then $\{g\in \mathbb F_k:\alpha_g\ne 0\}$. When viewing $x$ as an operator, we write  $x=\sum_{g\in \mathbb F_k}\alpha_g g$.
\begin{lemma}\label{FreeGroupLem}
With the notation above, let $H$ be a Hilbert space and $x\in M\ \vnotimes\ \mathbb B(H)$ satisfy $x(A\otimes \mathbb C1)x^*\subseteq A\ \vnotimes\ \mathbb B(H)$.  Then $x\in A\ \vnotimes\ \mathbb B(H)$.
\end{lemma}
\begin{proof}
With respect to some choice of matrix units for $\mathbb B(H)$, we may write $x = (x_{ij})_{i,j\in\Lambda}$ with $x_{ij} \in L(\mathbb F_k)$. The hypothesis implies that $\sum_jx_{ij}a^tx^*_{ji}\in A$ for all $t\in \mathbb Z$ and all $i\in \Lambda$, from which we wish to conclude that each $x_{ij}$ lies in $A$. Thus it suffices to consider operators $y_i\in L(\mathbb F_k)$ so that $\sum_jy_j a^ty^*_j\in A$ for all $t\in \mathbb Z$, and deduce that $y_j\in A$ for all $j$. Taking $t=0$ and applying the trace shows that $\sum_j \|y_j\|^2_2<\infty$, so by scaling we may assume that this sum is 1. We will argue by contradiction, so suppose that there is some   $j_0$ such that $y_{j_0}\notin A$. Then the support of $y_{j_0}$ contains a word $g$ which is not a power of $a$.   Then $g$ may be written in reduced form as $a^nw_0a^m$ where $m,n\in \mathbb Z$ and $w_0$ is a non-trivial reduced word whose first and last letters lie in $\{b_1^{\pm1},\dots,b_{k-1}^{\pm 1}\}$. We will examine the coefficient of $ga^tg^{-1}$ in $S(t):=\sum_jy_ja^ty_j^*$ and show that it is nonzero for sufficiently large values of $t$. This will give the desired contradiction.

Let $c=\langle y_{j_0}\delta_e,\delta_g\rangle\neq 0$ be the the coefficient of $g$ in $y_{j_0}$. For any fixed group element $h\in \mathbb F_k$ and any $t\in \mathbb Z$, the coefficient of $h$ in $y_ja^ty^*_j$ is
\begin{equation}
 \langle y_ja^ty^*_j\delta_e,\delta_h\rangle = \langle a^ty^*_j\delta_e, Jh^{-1}Jy^*_j\delta_e\rangle,
\end{equation}
where $Jh^{-1}J$ is the unitary operator of right-convolution by $h$.  This is bounded in absolute value by $\|y_j\|^2_2$ using the Cauchy-Schwarz inequality. Consequently, the contribution of $\sum_{j\in \Lambda_0} y_ja^ty^*_j$ to the coefficient of $ga^tg^{-1}$ in $S(t)$ is bounded above by $\sum_{j\in \Lambda_0}\|y_j\|^2_2$, for any subset $\Lambda_0$ of the index set $\Lambda$. Choose a finite set $\Lambda_0\subseteq\Lambda$ so that $\sum_{j\in\Lambda\setminus \Lambda_0}\|y_j\|^2_2 < |c|^2/4$. For $t\in\mathbb Z$, write $S_{\Lambda_0}(t)=\sum_{j\in \Lambda_0}y_ja^ty^*_j$.

For each $j\in \Lambda_0$, we may write $y_j$ as an orthogonal sum $z_j+z'_j$ where each $z_j$ is a finite linear combination of group elements, $g$ appears in $z_{j_0}$, and $\sum_{j\in \Lambda_0}\|z'_j\|_2 < |c|^2/4$. Then $ga^tg^{-1}$ can appear in $S_{\Lambda_0}(t)$ from the four terms $\sum_{j\in \Lambda_0}z_ja^tz^*_j$, $\sum_{j\in \Lambda_0}z_ja^tz^{\prime*}_j$, $\sum_{j\in \Lambda_0}z'_ja^tz^*_j$ and $\sum_{j\in \Lambda_0}z'_ja^tz^{\prime*}_j$. Since $\|z'_j\|_2 \le \|y_j\|_2\le 1$ and $\|z_j\|_2\leq\|y_j\|_2\leq 1$, we may argue as above to see that the total contribution (in absolute value) of the latter three terms to the coefficient of $ga^tg^{-1}$ is at most $3|c|^2/4$. For example, the contribution of the term $\sum_{j\in \Lambda_0}z_ja^tz^{\prime*}_j$ is bounded above in absolute value by
\begin{equation}
\sum_{j\in \Lambda_0} \|z_j\|_2\|z^{\prime*}_j\|_2\leq \sum_{j\in \Lambda_0} \|z^{\prime*}_j\|_2<|c|^2/4.
\end{equation}
We now examine the contribution from $\sum_{j\in \Lambda_0} z_ja^tz^*_j$, recalling that each $z_j$ is supported on finitely many group elements. Thus there is an integer $K>0$ which bounds the number of $a^{\pm 1}$'s  in any word in the support of $z_j$ for any $j\in \Lambda_0$.  There are two forms for such words. The first is $a^p$ where $|p|\leq K$ while the second is $a^p v a^q$, where $p,q\in \mathbb Z$ and the first and last letters of $v$ lie in $\{b_1^{\pm1},\dots,b_{k-1}^{\pm1}\}$. Moreover, we will have that $|p|,|q|\le K$, and that $v$ can contain at most $K$ $a^{\pm 1}$'s. We now restrict attention to $t\ge 5K$. In order to avoid a proliferation of cases, we write the typical term in the expansion of $\sum_{j\in \Lambda_0} z_ja^tz^*_j$ as $a^pva^qa^ta^{-s}w^{-1}a^{-r}$, where $|p|,|q|,|s|,|r|\le K$, each of $v$ and $w$ is either $e$ or begins and ends with letters from $\{b_1^{\pm1},\dots,b_{k-1}^{\pm1}\}$ and contains at most $K$ $a^{\pm1}$'s. We wish to know when this equals $a^nw_0a^ma^ta^{-m}w^{-1}_0a^{-n} = a^nw_0a^tw^{-1}_0a^{-n}$. This is impossible when both $v$ and $w$ are $e$, so we first consider the degenerate case when $v=e$ and the first and last letters of $w$ lie in $\{b_1^{\pm1},\dots,b_{k-1}^{\pm1}\}$. Then the equation
\begin{equation}\label{eq1}
a^{p+q+t-s}w^{-1}a^{-r}=a^nw_0a^tw^{-1}_0a^{-n}
\end{equation}
is false since $w_0$ starts with a letter from $\{b_1^{\pm1},\dots,b_{k-1}^{\pm1}\}$ while $p+q+t-s-n\geq K$. A similar argument disposes of the possibility that $w=e$. Thus we need only consider the case when both $v$ and $w$ begin and end with letters from $\{b_1^{\pm1},\dots,b_{k-1}^{\pm1}\}$, giving $p=r=n$. This question then reduces to examining the equation
\begin{equation}\label{eq2}
 va^{t+q-s} w^{-1} = w_0a^tw^{-1}_0.
\end{equation}
The last letter in $v$ either lies in $w_0$ or in $w^{-1}_0$. In the second case $v$ must contain at least $w_0a^t$, contradicting the bound of $K$ on the number of $a$'s in $v$ since $t\geq 5K$. Thus the first case must hold, and we can write $w_0=vw_1$ where no cancelations can occur between $v$ and $w_1$. Then (\ref{eq2}) becomes
\begin{equation}\label{eq3}
 a^{t+q-s}w^{-1} = w_1a^tw^{-1}_0.
\end{equation}
Since the last letter of $w_0$ lies in $\{b_1^{\pm1},\dots,b_{k-1}^{\pm1}\}$, there is no cancelation between $a^t$ and $w_0^{-1}$. Thus there are two possibilities:\ $w_1=e$ or $w_1$ ends with a letter from $\{b_1^{\pm1},\dots,b_{k-1}^{\pm1}\}$. If the second one holds then no cancelations occur on the right hand side of (\ref{eq3}), so $w_1$ begins with $a^{t+q-s}$, showing that $w_0$ contains this power. But $t+q-s\ge 3K$ and a contradiction ensues. Thus $w_1=e$ and $w_0=v$. Returning to (\ref{eq2}), we conclude that $q=s$ and $w_0=w$.  It follows that the contributions to the coefficient of $ga^tg^{-1}$ in $\sum_{j\in \Lambda_0}z_ja^tz_j^*$ only arise from products $da^td^{-1}$ where $d$ has the form $a^nw_0a^q$. Thus this coefficient is a sum of positive terms including $|c|^2$ (which arises from the coefficient $ca^nw_0a^m$ term in $z_{\Lambda_0}$).  In conclusion, for $t$ sufficiently large, the coefficient of $ga^tg^{-1}$ in $\sum_{j\in \Lambda} y_ja^ty_j^*$ has absolute value at least $|c|^2-|c|^2/4-3|c|^2/4$, and is thus nonzero.  Consequently these elements cannot lie in $A$, and this contradiction completes the proof.
\end{proof}

\begin{theorem}\label{FreeGroup}
Let $k\geq 2$ and let $A$ be a generator masa in $M=L\mathbb F_k$.  Given any  unital inclusion $B\subseteq N$ of von Neumann algebras, we have
\begin{equation}
\GN_{M\,\vnotimes\,N}(A\ \vnotimes\ B)''=A\ \vnotimes\ \GN_{N}(B)''.
\end{equation}
\end{theorem}
\begin{proof}
Take a faithful representation of $N$ on some Hilbert space $H$.  Let $v$ be a groupoid normaliser of $A\ \vnotimes\ B$ in $M\ \vnotimes\ N$. Then $v(A\otimes\mathbb C1)v^*\subseteq A\ \vnotimes\ B\subseteq A\ \vnotimes\ \mathbb B(H)$, so by Lemma \ref{FreeGroupLem}, $v$ lies in $A\ \vnotimes\ \mathbb B(H)$.  By applying slice maps, it follows that $v\in A\ \vnotimes\ N$. Thus
\begin{equation}
\GN_{M\,\vnotimes\,N}(A\ \vnotimes\ B)''=\GN_{A\,\vnotimes\,N}(A\ \vnotimes\ B)'',
\end{equation}
as the inclusion from right to left is immediate. Lemma \ref{AbDirectInt} then gives
\begin{equation}
\GN_{M\,\vnotimes\,N}(A\ \vnotimes\ B)''=A\ \vnotimes\ \GN_N(B)'',
\end{equation}
and so the result follows.
\end{proof}

\section{Irreducible subfactors}\label{Irred}

Recall that an inclusion $B\subseteq M$ of factors is \emph{irreducible} if $B'\cap M=\mathbb C1$.  In this section we consider two irreducible inclusions $B_i\subseteq M_i$ ($i=1,2$) of factors and  the resulting tensor product inclusion $B=B_1\ \vnotimes\ B_2\subseteq M_1\ \vnotimes\ M_2=M$.  Our objective  is to examine normalisers of this  latter inclusion and to address the question of whether  they always factorise in the form $v(w_1\otimes w_2)$ for some normalisers $w_i\in\N_{M_i}(B_i)$ and some unitary $v\in B$. In \cite{Saw.NormalisersSubfactors} a positive answer  was obtained when both $M_1$ and $M_2$ are type II$_1$ using basic construction methods (the special case of finite index inclusions follows from the earlier paper \cite{Popa.Entropy}). Here we examine the $M'$-bimodular automorphism groups used in Sections \ref{Aut} and \ref{Tensor} to answer this question beyond the finite setting.  To this end, we regard each $M_i$ as faithfully represented on a Hilbert space $H_i$ so that $M$ is represented on $H=H_1\otimes H_2$. Tomita's commutation theorem ensures that the commutants $M'$ and $B'$ are obtained as the tensor products $M_1'\ \vnotimes\ M_2'$ and $B_1'\ \vnotimes\ B_2'$ respectively.  Given any unitary normaliser $u\in\N_M(B)$, we obtain an automorphism $\theta=\ad(u)\in\Aut_{M'}(B')$.  In the setting of irreducible subfactors, these automorphisms must factorise, as we now show.

\begin{lemma}\label{Irred.Split}
With the notation  above, let $\theta\in\Aut_{M'}(B')$. Then there exist automorphisms $\theta_i\in\Aut_{M_i'}(B_i')$, $i=1,2$, such that $\theta=\theta_1\otimes\theta_2$.
\end{lemma}
\begin{proof}
For $y_1\in B_1'$ and $x_2\in M_2'$, we have
\begin{align}
(1\otimes x_2)\theta(y_1\otimes 1)&=\theta((1\otimes x_2)(y_1\otimes 1))=\theta((y_1\otimes 1)(1\otimes x_2))\notag\\
&=\theta(y_1\otimes 1)(1\otimes x_2)
\end{align}
so that 
\begin{equation}
\theta(y_1\otimes 1)\in (B_1'\ \vnotimes\ B_2')\cap (\mathbb C1\ \vnotimes M_2')'=B_1'\ \vnotimes \ (B_2'\cap M_2)=B_1'\otimes\mathbb C1.
\end{equation}
Similarly, $\theta^{-1}$ also maps $B_1'\otimes \mathbb C1$ into itself, so if we identify $B_1'$ with $B_1'\otimes\mathbb C1$, then $\theta$ restricts to give an automorphism of $B_1'$.  By interchanging the roles of $B_1$ and $B_2$, we see that  $\theta$ also restricts to an automorphism $\theta_2$ of $B_2'$, allowing us to conclude that $\theta=\theta_1\otimes\theta_2$.  Given $x_1\in M_1'$ and $x_2\in M_2'$, we have 
\begin{equation}
x_1\otimes x_2=\theta(x_1\otimes x_2)=\theta_1(x_1)\otimes \theta_2(x_2),
\end{equation}
so that each $\theta_i$ is $M_i'$-bimodular for $i=1,2$.
\end{proof}

The next lemma enables us to recover groupoid normalisers from elements commuting with bimodular automorphisms.
\begin{lemma}\label{Irred.Reg}
Let $B\subseteq M$ be an irreducible inclusion of subfactors, faithfully represented on some Hilbert space $H$ and let $\psi$ be an element of $\Aut_{M'}(B')$.  Let $a\in\mathbb B(H)$ have polar decomposition $a=v|a|$ and also satisfy
\begin{equation}\label{Irred.Reg.1}
\psi(y)a=ay,\quad y\in B'.
\end{equation}
 Then $aa^*,a^*a\in B$ and $v$ is a groupoid normaliser of $B$ in $M$ satisfying
\begin{equation}
\psi(y)v=vy,\quad y\in B'.
\end{equation}
\end{lemma}
\begin{proof}
Applying (\ref{Irred.Reg.1}) and its adjoint gives
\begin{equation}
a^*ay=a^*\psi(y)a=ya^*a,\quad y\in B',
\end{equation}
so that $a^*a\in B$. Similarly
\begin{equation}
aa^*\psi(y)=\psi(y)aa^*,\quad y\in B'
\end{equation}
so that $aa^*\in B$. Now,  
\begin{equation}
\psi(y)v|a|=v|a|y=vy|a|,\quad y\in B',
\end{equation}
which implies that
\begin{equation}\label{Irred.Reg.2}
\psi(y)v=vyv^*v=vv^*vy=vy,\quad y\in B'
\end{equation}
since $v^*v$ is the range projection of $|a|$ (and so lies in $B$).  By hypothesis $\psi(x)=x$ for $x\in M'$, so the double commutant theorem ensures that $v\in M$.  For $b\in B$ and $y\in B'$, we have
\begin{equation}
vbv^*\psi(y)=vbyv^*=vybv^*=\psi(y)vbv^*,
\end{equation}
so that $vbv^*\in B$.  Similarly
\begin{equation}
v^*bvy=v^*b\psi(y)v=v^*\psi(y)bv=yv^*bv,
\end{equation}
so $v^*bv\in B$ and we have proved that $v\in\GN_M(B)$.
\end{proof}

We briefly recall from \cite{Tom.Cop} some facts about slice maps on tensor
products that will be used subsequently. If $P$ and $Q$ are two von Neumann
algebras then the algebraic tensor product $P_*\otimes Q_*$ is a norm dense
subspace of the predual of $P\overline{\otimes }Q$. Thus any fixed $\phi\in
Q_*$ defines a left slice map $L_\phi
:P\ \overline{\otimes}\ Q\to P$ by
\begin{equation}
L_\phi(x)(\psi)=(\psi\otimes\phi)(x),\ \ \ x\in P\ \overline{\otimes}\ Q,\ \ \ 
\psi\in P_*.
\end{equation}
This definition shows that each $L_\phi$ is a normal map and, when
restricted to elementary tensors, it has the form
\begin{equation}\label{elem}
L_\phi(p\otimes q)=\phi(q)p,
\ \ \  p\in P,\ q\in Q. \end{equation}
It is immediate from (\ref{elem}) that each $L_\phi$ is a $P$-bimodule map.
Moreover, if $N$ is a von Neumann subalgebra of $P$ and $x\in
P\ \overline{\otimes}\ Q$ is such that $L_\phi(x)\in N$ for all $\phi \in Q_*$,
then $x\in N\ \overline{\otimes}\ Q$. In the next lemma we note that the automorphism $\theta$ splits as $\theta=\theta_1\otimes \theta_2$, by Lemma \ref{Irred.Split}.

\begin{lemma}\label{Irred.Tech}
For $i=1,2$, let $B_i\subseteq M_i$ be irreducible inclusions of factors. Write $M=M_1\ \vnotimes\ M_2$ and $B=B_1\ \vnotimes\ B_2$ and let $u\in \N_{M}(B)$.  Write $\theta=\theta_1\otimes \theta_2$ for the induced element $\ad(u)$ of $\Aut_{M'}(B')$. Then there exist $v_i\in\GN_{M_i}(B_i)$,  $i=1,2$, such that each $v_i$ is either an isometry or a co-isometry and 
\begin{equation}
\theta_i(y_i)v_i=v_iy_i,\quad y\in B_i',\ i=1,2.
\end{equation}
\end{lemma}
\begin{proof}
By symmetry, we only need examine the case $i=1$. Consider the set 
\begin{equation}
S_1=\{v\in\GN_{M_1}(B_1):\theta_1(y)v=vy,\quad y\in B_1'\},
\end{equation}
which we equip with the partial ordering $v\leq w$ iff $wv^*v=v$ (equivalent to the requirement that  the partial isometry $w$ be an extension of $v$).  It is clear that any chain in $S_1$ has a least upper bound so, by Zorn's Lemma, there is some maximal element $v\in S_1$.  Let $e=1-v^*v$ and $f=1-vv^*$.  Suppose that $v$ is neither an isometry nor a co-isometry, so that both $e$ and $f$ are non-zero.  

Since $u(e\otimes 1)\neq 0$, there is a normal state $\phi$ on $\mathbb B(H_2)$  
  inducing a slice map $\mathrm{id}\otimes \phi:\mathbb B(H_1)\ \vnotimes\ \mathbb B(H_2)\rightarrow\mathbb B(H_1)$ with $(\mathrm{id}\otimes\phi)(u(e\otimes 1))\neq 0$.   If we define $a\in \mathbb B(H_1)$ by $a=(\mathrm{id}\otimes \phi)(u(e\otimes 1))$, then the module property of the slice map shows that 
\begin{equation}
a=(\mathrm{id}\otimes \phi)(u(e\otimes 1))=(\mathrm{id}\otimes \phi)(u))e.
\end{equation}
  Now
\begin{equation}
(\theta_1(y)\otimes 1)u=\theta(y\otimes 1)u=u(y\otimes 1),\quad y\in B_1',
\end{equation}
so applying $\mathrm{id}\otimes\phi$ gives
\begin{equation}
\theta_1(y)(\mathrm{id}\otimes \phi)(u)=(\mathrm{id}\otimes\phi)(u)y,\quad y\in B_1'.
\end{equation}
Then multiplication on the right by $e\in B_1$ leads to 
\begin{equation}
\theta_1(y)a=ay,\quad y\in B_1'.
\end{equation}
Letting $a=v_1|a|$ be the polar decomposition of $a$, we  obtain $v_1\in S_1$ from Lemma \ref{Irred.Reg}. Define $f_1=v_1v_1^*$, and note that    some non-zero subprojection of $f_1$ is equivalent to a subprojection of $f$ since $B_1$ is a factor.  Let $w\in B_1$ be a nonzero partial isometry with $w^*w\leq f_1$ and $ww^*\leq f$, and consider the partial isometry $wv_1$.  This is a groupoid normaliser of $B_1$ with
\begin{equation}
\theta_1(y)wv_1=wv_1y,\quad y\in B_1',
\end{equation}
since both $v_1$ and $w$ have these properties.  Accordingly, $v+wv_1$ lies in $S_1$ and is strictly greater than $v$,  contradicting maximality.  Hence $v$ must be either an isometry or a co-isometry.
\end{proof}

In the following corollary we recapture the normaliser results for tensor products of irreducible inclusions of II$_1$ factors from \cite{Saw.NormalisersSubfactors}.
\begin{corollary}\label{Irred.Cor}
For $i=1,2$, let $B_i\subseteq M_i$ be irreducible inclusions of ${\mathrm{II}}_1$ factors. Then any unitary normaliser of $B_1\ \vnotimes\ B_2$ in $M_1\ \vnotimes\ M_2$ is of the form $w(v_1\otimes v_2)$ for some unitary $w\in B_1\ \vnotimes\ B_2$ and normalisers $v_i\in\N_{M_i}(B_i)$. 
\end{corollary}
\begin{proof}
Given a normaliser $u$ of $B_1\ \vnotimes\ B_2$ in $M_1\ \vnotimes\ M_2$, let $\theta=\ad(u)\in\Aut_{M_1'\,\vnotimes\,M_2'}(B_1'\ \vnotimes\ B_2')$ (formed with respect to some faithful actions of $M_1$ and $M_2$ on Hilbert spaces $H_1$ and $H_2$). By Lemma \ref{Irred.Split}, we have a factorisation $\theta=\theta_1\otimes\theta_2$ for some $\theta_i\in\Aut_{M_i'}(B_i')$.  Lemma \ref{Irred.Tech} gives groupoid normalisers $v_i\in\GN_{M_i}(B_i)$ satisfying
\begin{equation}
\theta_i(y_i)v_i=y_iv_i,\quad y_i\in B_i',
\end{equation}
and such that $v_1$ and $v_2$ are either isometries or co-isometries.  Since both $M_1$ and $M_2$ are finite, $v_1$ and $v_2$ are unitary normalisers.  Define $w=u(v_1^*\otimes v_2^*)$.  Since
\begin{equation}
\theta(y_1\otimes y_2)w=u(y_1v_1^*\otimes y_2v_2^*)=u(v_1^*\theta_1(y_1)\otimes v_2^*\theta_2(y_2))=w\theta(y_1\otimes y_2),\quad y_i\in B_i',
\end{equation}
the double commutant theorem ensures that $w\in B_1\ \vnotimes\ B_2$ and $u=w(v_1\otimes v_2)$, as required.
\end{proof}

We now return to the situation where we do not assume that $M_1$ and $M_2$ are finite.  
\begin{lemma}\label{Lem5.5}
Let $B\subseteq M$ be an irreducible inclusion of infinite factors and suppose that $\psi\in\Aut_{M'}(B')$ is given and $v\in\GN_M(B)$ is an isometry or co-isometry with 
\begin{equation}\label{Irred.Infinite.2}
\psi(y)v=vy,\quad y\in B'.
\end{equation}
Then there exists a normaliser $v_1\in\N_M(B)$ with 
\begin{equation}\label{Irred.Infinite.1}
\psi(y)v_1=v_1y,\quad y\in B'.
\end{equation}
\end{lemma} 
\begin{proof}
We may assume that $v^*v=1$ and $vv^*=f<1$. The map $x\mapsto vxv^*$ implements a surjective $^*$-isomorphism between $B$ and $fBf$. Hence $fBf$ is an infinite factor so $f$ is infinite in $B$. Then there exists a partial isometry $w\in B$ with $w^*w=f$ and $ww^*=1$. Define $v_1=wv$. This is certainly a unitary normaliser of $B$ and (\ref{Irred.Infinite.1}) follows from multiplying (\ref{Irred.Infinite.2}) on the left by $w$ (which commutes with $\psi(y)$).
\end{proof}

\begin{theorem}\label{IrredThm1}
For $i=1,2$, let $B_i\subseteq M_i$ be irreducible inclusions of factors.  Suppose that either one of $M_1$ or $M_2$ is finite or that both $B_1$ and $B_2$ are infinite.  Then any unitary normaliser $u$ of $B_1\ \vnotimes\ B_2$ in $M_1\ \vnotimes\ M_2$ is of the form $w(v_1\otimes v_2)$ for some unitary $w\in B_1\ \vnotimes\ B_2$ and normalisers $v_i\in\N_{M_i}(B_i)$. 
\end{theorem}
\begin{proof}
As in the proof of Corollary \ref{Irred.Cor}, write $\ad(u)=\theta_1\otimes\theta_2$ for some $\theta_i\in\Aut_{M_i'}(B_i')$, and find isometries or co-isometries $v_i\in\GN_{M_i}(B_i)$ with 
\begin{equation}
\theta_i(y_i)v_i=y_iv_i,\quad y_i\in B_i'.
\end{equation}
If both $B_1$ and $B_2$ are infinite, then Lemma \ref{Lem5.5} enables us to replace the $v_i$ by unitary normalisers satisfying the same condition. Just as in Corollary \ref{Irred.Cor}, $w=u(v_1^*\otimes v_2^*)$ is a unitary in $B_1\ \vnotimes\ B_2$ and in this case $u$ has the desired form.  Now consider the case where   $M_1$ is finite. Then $v_1$ is automatically a unitary normaliser. If $B_2$ is infinite, then we can use Lemma \ref{Lem5.5} to replace $v_2$ by a normaliser, and the result follows. If $B_2$ is finite, then the element $w$, defined by $w=u(v_1^*\otimes v_2^*)$,  lies in $B_1\ \vnotimes\ B_2$ and is either an isometry or co-isometry (since $u$ is a unitary).  In this case, the tensor product $B_1\ \vnotimes\ B_2$ is finite, so $w$, and hence $v_2$, is a unitary as required.
\end{proof}

To consider the last remaining case when both inclusions $B_i\subseteq M_i$ consist of an irreducible inclusion of a II$_1$ factor inside an infinite factor, we use Murray and von Neumann's fundamental group from \cite{MvN.4}. Let $B$ be a II$_1$ factor with trace $\tau$.  Recall that the fundamental group, $\mathcal F(B)$, of $B$ is the subgroup of $\mathbb R^+$ of those quotients $\tau(p)/\tau(q)$, where $p$ and $q$ are projections in $B$ such that $pBp$ and $qBq$ are isomorphic. Many subgroups of $\mathbb{R}^+$ can occur, in particular all countable subgroups of $\mathbb R^+$ \cite{Popa.StrongRigidity1}.  In our situation, it is enough for one subfactor $B_i$ to have trivial fundamental group; this enables us to convert the isometries or co-isometries appearing in the proof of Theorem \ref{IrredThm1} into unitary normalisers.
\begin{lemma}\label{FG.1=>3}
For $i=1,2$, let $B_i\subseteq M_i$ be irreducible inclusions of factors. Suppose that $B_1$ is finite and has trivial fundamental group. Then any unitary normaliser of $B_1\ \vnotimes\ B_2$ in $M_1\ \vnotimes\ M_2$ is of the form $w(v_1\otimes v_2)$ for some unitary $w\in B_1\ \vnotimes\ B_2$ and normalisers $v_i\in\N_{M_i}(B_i)$. 
\end{lemma}
\begin{proof}
As in the proof of Corollary \ref{Irred.Cor}, write $\ad(u)=\theta_1\otimes\theta_2$ for some $\theta_i\in\Aut_{M_i'}(B_i')$, and find isometries or co-isometries $v_i\in\GN_{M_i}(B_i)$ with 
\begin{equation}
\theta_i(y_i)v_i=y_iv_i,\quad y_i\in B_i'.
\end{equation}
If $v_1$ is an isometry, then $x\mapsto v_1xv_1^*$ implements an isomorphism of $B$ onto $(v_1v_1^*)B(v_1v_1^*)$. As $\mathcal F(B_1)=\{1\}$, we have $v_1v_1^*=1$ so $v_1$ is a unitary. If $v_1$ is a co-isometry, we can consider $x\mapsto v_1^*xv_1$ to see that $v_1$ is again a unitary.  The rest of the argument now follows  exactly the proof of Theorem \ref{IrredThm1}, considering the cases where $B_2$ is infinite and $B_2$ is finite separately.
\end{proof}

In Theorem \ref{thm5.9} we will establish a converse to Lemma \ref{FG.1=>3}.
First fix a finite factor $B$ and then consider arbitrary irreducible inclusions $B\subseteq M$. We will show that if the conclusion of Lemma \ref{FG.1=>3}
holds for the tensor product with every other separable irreducible
inclusion $B_2\subseteq M_2$, then $B$ has trivial fundamental group.

Recall that if $B$ is a II$_1$ factor with trace $\tau_B$, then $\langle
x,y\rangle=\tau_B(y^*x)$ defines an inner product on $B$, inducing a norm
$\|x\|_2=\tau_B(x^*x)^{1/2}$ for $x\in B$. The completion is denoted
$L^2(B)$, we write $\xi$ for the image of 1 in this Hilbert space, and
$\{x\xi:x\in B\}$ is a dense subspace. The representation of $B$ by
$b(x\xi)=(bx)\xi$ for $b, x\in B$, extends by continuity to a representation
of $B$ on $L^2(B)$ and puts $B$ into standard form on this Hilbert space.

Now suppose that $B$ is a II$_1$ factor with non-trivial fundamental group. Fix a projection $e\in B$ with $e\neq 0,1$ such that $B\cong eBe$ and choose an isomorphism $\theta$ from $B$ onto $eBe$.  The trace on $eBe$ is given by $\tau_{eBe}(x)=\tau_B(x)/\tau_B(e)$ and since this trace is unique, we must have
\begin{equation}
\tau(x)=\tau_{eBe}(\theta(x)),\quad x\in B.
\end{equation}
Thus 
\begin{equation}
\|\tau_B(e)^{-1/2}\theta(x)\|_2^2=\frac{\tau_B(\theta(x^*x))}{\tau_B(e)}=\tau_{eBe}(\theta(x^*x))=\tau_B(x^*x)=\|x\|_2^2,\quad x\in B.
\end{equation}
Therefore, we can define an isometry $v$ on $L^2(B)$, by extending
\begin{equation}
v(x\xi)=\frac{1}{\tau(e)^{1/2}}\theta(x)\xi,\quad x\in B,
\end{equation}
by continuity.  This operator implements $\theta$, in that $\theta(x)=vxv^*$ for all $x\in B$ and as $\theta$ is a surjective isomorphism, $ebe\mapsto v^*bv$ must be the inverse of $\theta$. Since $vv^*=e$, this shows that $v^*Bv=B$ so $v$ is a groupoid normaliser of $B$ in $\mathbb B(L^2(B))$. Now let $H=L^2(B)\otimes\ell^2(\mathbb Z)$.  Let $v_1$ be the bilateral shift operator on $\ell^2(\mathbb Z)$ satisfying $v_1e_n=e_{n+1}$ for  $n\in\mathbb Z$ and define $u=v\otimes v_1\in\mathbb B(H)$.  Let $M$ be the von Neumann subalgebra of $\mathbb B(H)$ generated by $B\otimes 1$ and $u$.
This algebra resembles a  crossed product; indeed it would be $B\rtimes_{\alpha}\mathbb{Z}$ for the action induced by $\ad(u)$ if $u$ were a unitary. 

\begin{theorem}\label{thm5.8}
In the situation described above, $B\otimes\mathbb C1\subseteq M$ is a singular inclusion of factors with $\GN_M(B\otimes 1)''=M$.
\end{theorem}
\begin{proof}
By construction $v$ is a groupoid normaliser of $B$ with $v(B')v^*=B'e$. Since $B'$ is isomorphic to $B'e$ via $b'\mapsto b'e$, we can define an automorphism $\phi$ of $B'$ so that $\phi(b')e=vb'v^*$.  Given some $b'\in B'$, define a bounded operator $y=(y_{i,j})_{i,j}$ on $H$ with respect to the matrix units for $\ell^2(\mathbb Z)$ by setting $y_{i,i}=\phi^i(b')$ for all $i$ and setting the off-diagonal elements to be zero. Certainly $y$ commutes with $(B\otimes\mathbb C1)$.  Now, for $n\in\mathbb Z$ and $\eta\in L^2(B)$,
\begin{align}
yu(\eta\otimes e_n)&=y(v\eta \otimes e_{n+1})
=\phi^{n+1}(b')ev\eta\otimes e_{n+1}\notag\\
&=v\phi^n(b')v^*v\eta\otimes e_{n+1}
=v\phi^n(b')\eta\otimes e_{n+1}\notag\\
&=uy(\eta\otimes e_n),
\end{align}
so $y$ commutes with  $u$. The elements $y$ that we are considering form a self-adjoint set, and so they also commute with $u^*$. Thus these elements  commute with the generators of $M$ and so lie in $M'$.

Now take $x\in M$, let $x=(x_{i,j})$ with respect to the matrix units of $\ell^2(\mathbb Z)$.  Since $1\otimes v_1$ commutes with $B\otimes \mathbb C1$ and $u=v\otimes v_1$, $1\otimes v_1$ lies in $M'$.  Thus $x=(1\otimes v_1)^*x(1\otimes v_1)$. This shows that $x$ has constant diagonal entries, i.e. $x_{i,i}=x_{0,0}$ for all $i$. Given $b'\in B'$, define $y\in M'$ as above.  Comparing the $(0,0)$-th entries of $xy=yx$, we see that $x_{0,0}b'=b'x_{0,0}$ so that $x_{0,0}=b_0$ for some $b_0\in B$. For $n\geq 0$, consider $x(v^*\otimes v_1^*)^n$ which lies in $M$.  Thus this operator has constant entries, say $b_{-n}\in B$, down the diagonal and so $x_{i+n,i}=b_{-n}v^n$ for all $i$.  Similarly we can find bounded operators $b_n\in B$ for $n>0$ such that $x_{i-n,i}=(v^*)^nb_n$.  Then the general form of an operator $x\in M$ is
\begin{equation}\label{101.1}
x=\sum_{n=0}^\infty (b_{-n}v^n\otimes v_1^n)+\sum_{n=1}^\infty ((v^*)^nb_n\otimes (v_1^*)^n),
\end{equation}
with convergence in the weak operator topology.  For $n\geq 1$, define $e_n=\theta(e_n)=v^n(v^*)^n$. Since $b_{-n}v^n=b_{-n}e_nv^n$, we may assume that $b_{-n}=b_{-n}e_n$ in (\ref{101.1}) for $n>0$. Similarly, we can assume that $b_n=e_nb_n$ for $n>0$.

Now suppose that $x$ is a normaliser of $B\otimes \mathbb C1$ in $M$ and write $x$ in the form (\ref{101.1}).  Define $\psi:B\rightarrow B$ by $\psi(b)\otimes 1=x^*(b\otimes 1)x$ so that $(b\otimes 1)x=x(\psi(b)\otimes 1)$ for $b\in B$. Substituting  this into (\ref{101.1}) gives
\begin{equation}\label{101.1a}
\sum_{n=0}^\infty (bb_{-n}v^n-b_{-n}v^n\psi(b))\otimes v_1^n+\sum_{n=1}^\infty (b(v^*)^nb_n-(v^*)^nb_n\psi(b))\otimes (v_1^*)^n)=0
\end{equation}
for each $b\in B$.  Comparing matrix elements in (\ref{101.1a}) leads to
\begin{equation}\label{101.2}
bb_{-n}v^n=b_{-n}v^n\psi(b),\quad n\geq 0,\ b\in B,
\end{equation}
and
\begin{equation}\label{101.3}
b(v^*)^nb_n=(v^*)^nb_n\psi(b),\quad n>1,\ b\in B.
\end{equation}
After multiplying (\ref{101.2}) on the right by $(v^*)^n$, we see that
\begin{equation}
bb_{-n}=b_{-n}\theta^n(\psi(b)),\quad n\geq 0,\ b\in B.  
\end{equation}
When $b=b^*$, we can apply this twice to obtain
\begin{equation}
bb_{-n}b_{-n}^*=b_{-n}\theta(\psi(b))b_{-n}^*=b_{-n}b_{-n}^*b,\quad n\geq 0,\ b=b^*\in B,
\end{equation}
so that $b_{-n}b_{-n}^*\in B'\cap B=\mathbb C1$. Thus $b_{-n}$ is a scalar multiple of a unitary in $B$. Then the relations $b_{-n}=b_{-n}e_n$ and $e_n\neq 1$ for $n>0$ imply that $b_{-n}=0$ for $n>0$.  Replacing $x$ by $x^*$ shows that $b_n=0$ for $n<0$ and so $x=b_0\otimes 1\in B\otimes\mathbb C1$. Thus $B\otimes\mathbb C 1$ is singular in $M$, so
$(B\otimes\mathbb C 1)'\cap M\subseteq (B\otimes \mathbb C1)'\cap B\otimes
\mathbb C 1=\mathbb C1$. In particular, $M$ has trivial centre and so we have an inclusion of factors. Finally, $M$ is generated by $B\otimes \mathbb C1$ and the groupoid normaliser $u=v\otimes v_1$ so that $M=\GN_{M}(B\otimes \mathbb C1)''$.  
\end{proof}
Our final result characterises the property of having trivial fundamental group in terms of normalisers of tensor products.

\begin{corollary}\label{thm5.9}
Let $B$ be a ${\mathrm{II}}_1$ factor with a separable predual.  The following statements are equivalent.
\begin{enumerate}[(i)]
\item The fundamental group of $B$ is trivial.\label{5.9L1}
\item Whenever $B\subseteq M$ is an inclusion of factors, $\GN_M(B)''=\N_M(B)''$.\label{5.9L2}
\item Whenever $B\subseteq M$ is an irreducible inclusion of factors and $B_2\subseteq M_2$ is another irreducible inclusion of factors, then every unitary normaliser $u\in\N_{M\,\vnotimes\,M_2}(B\ \vnotimes\ B_2)$ factorises as $u=w(u_1\otimes u_2)$, where $w\in B\ \vnotimes\ B_2$, $u_1\in\N_M(B)$ and $u_2\in\N_{M_2}(B_2)$.\label{5.9L3}
\end{enumerate}
\end{corollary}
\begin{proof}
Lemma \ref{FG.1=>3} shows that $(\ref{5.9L1})\Rightarrow (\ref{5.9L3})$ and the previous theorem establishes $(\ref{5.9L2})\Rightarrow (\ref{5.9L1})$, so it remains to show $(\ref{5.9L3})\Rightarrow (\ref{5.9L2})$.  Take $B_2=M_2=\mathbb B(\ell^2(\mathbb N))$ so that condition $(\ref{5.9L3})$ implies that
\begin{equation}
\N_{M\,\vnotimes\,\mathbb B(\ell^2(\mathbb N))}(B\ \vnotimes\ \mathbb B(\ell^2(\mathbb N)))=\N_M(B)''\ \vnotimes\ \mathbb B(\ell^2(\mathbb N)).
\end{equation}
Thus $\GN_M(B)''=\N_M(B)''$ by Theorem \ref{GN=N}.
\end{proof}

\end{document}